\documentclass[leqno]{amsart}
\usepackage{amsmath,amsfonts,amssymb,amscd,amsthm,amsbsy,upref}
\usepackage{eufrak}
\newtheorem{thm}{Theorem}		
\newtheorem*{Thm}{Theorem}		
\newtheorem{lem}[thm]{Lemma}
\newtheorem{cor}[thm]{Corollary}
\newtheorem*{corA}{Corollary A}
\newtheorem*{corB}{Corollary B}
\newtheorem*{corC}{Corollary C}
\newtheorem*{propA}{Proposition A}
\newtheorem*{propB}{Proposition B}
\newtheorem{prop}[thm]{Proposition}
\newtheorem{prob}{Problem}
\theoremstyle{definition}

\newtheorem*{Define}{Definition}

\newtheorem*{Rem}{Remark}
\newtheorem*{remarks}{Remarks}
\def\ds{\displaystyle}
\def\bone{{\mathbf 1}}
\def\ep{\varepsilon}
\def\To{\Rightarrow}

\def\Ba{\operatorname{Ba}}
\def\co{\operatorname{co}}
\def\card{\operatorname{card}}
\def\dens{\operatorname{dens}}

\def\Ext{\operatorname{Ext}}
\def\Im{\operatorname{Im}}
\def\Re{\operatorname{Re}}
\def\tr{\operatorname{tr}}
\def\trinorm#1{{|\!|\!| #1 |\!|\!|}}
\def\vtimes{{\mathop{\otimes}\limits^\vee}}
\def\Atimes{{\mathop{\otimes}\limits^\wedge}}
\def\defeq{\ {\mathop{=}\limits^{\text{\rm def}}}\ }
\def\complex{{\mathbb C}}
\def\KK{{\mathbb K}}
\def\nat{{\mathbb N}}
\def\real{{\mathbb R}}
\def\B{{\mathcal B}}
\def\D{{\mathcal D}}
	\def\E{{\varepsilon}}
\def\F{{\mathcal F}}

\def\L{{\mathcal L}}
\def\P{{\mathcal P}}
\def\S{{\mathcal S}}
\def\U{{\mathcal U}}

\def\G{{\EuFrak c}}
\begin{document}
\title{Some new characterizations of Banach spaces containing $\ell^1$}
\author[Haskell P. Rosenthal]{Haskell P. Rosenthal*}
\begin{abstract} 
Several new characterizations of Banach spaces containing a subspace 
isomorphic to $\ell^1$, are obtained.
These are applied to the question of when $\ell^1$ embeds in the 
injective tensor product of two Banach spaces.
\end{abstract}
\thanks{*Research partially supported by NSF DMS-0700126.}

\maketitle

\subsection*{Notations and terminology}
All Banach spaces are taken as infinite dimensional, 
``subspace'' means ``closed linear subspace,'' 
``operator'' means ``bounded linear operator.''
If $W$ is a subset of a Banach space, $[W]$ denotes its closed linear span.
$\G$ denotes the cardinal of the continuum, i.e., $\G = 2^{\aleph_0}$; 
this is also identified with the least ordinal of cardinality $\G$.
For $1\le p<\infty$, $\ell_\G^p$ denotes the family $f$ of all scalar 
valued functions defined on $\G$ with 
$\|f\|_p = (\sum_{\alpha <\G} |f(\alpha)|^p)^{1/p} <\infty$. 
Finally, we recall that a scalar-valued function defined on a compact metric 
space $K$ is called universally measurable if it is measurable with respect 
to the completion of every Borel measure on $K$. 

Throughout this paper, the symbols $X,Y,Z,B,E$ shall denote Banach spaces. 
$\Ba X$ denotes the closed unit ball of $X$.
Recall that an operator $T:X\to Y$ is called Dunford-Pettis if $T$ maps 
weakly compact sets in $X$ to norm compact sets in $Y$. 
Also, $\L(X,Y)$ (resp. $K(X,Y)$) denotes the  space of operators
(resp. of compact operators) from $X$ to $Y$, 
$\L(X) = \L(X,X)$, $K(X) = K(X,X)$. 
A bounded subset $W$ of $X^*$ is said to {\em isomorphically norm\/} 
$X$ if there exists a $C< \infty$ such that
$$\|x\| \le C\sup_{w\in W} |w(x)| \ \text{ for all }\ x\in X\ .$$
In case $C= 1$ and $W\subset \Ba X^*$,   
we say that $W$ {\em isometrically norms} $X$. 
$X \mathop{\otimes}\limits^\vee Y$, $X \mathop{\otimes}\limits^\wedge Y$  
denote the injective, respectively projective, tensor products of $X$ and $Y$.
See \cite{DU}, \cite{Gr2} for terminology and theorems in this area.

\subsection*{Main results}
Our first main result gives several equivalences for a Banach space to 
contain an isomorph of $\ell^1$.
We have included many previously known ones, to round out the list; also,
we use some of them later on.
As far as I know, the equivalences of 1.\ with the following are new:
2, 3, 4, 6, 7, 8, 11, 12, 13, 14, 19.
(Of course some of the implications were previously known or are obvious. 
Also, the same construction proves that 2, 3, 4 and 6 imply 1, but I 
thought it useful nevertheless to list these explicitly.)
For other equivalences, cf.\ \cite{H1}, \cite{Hay}, \cite{G1}, 
Theorem~II.3 of \cite{G2}, \cite{Ro4}, and several of the remarks 
following the proof of Theorem~\ref{thm:main}. 

\begin{thm}\label{thm:main}
Let $X$ be given. 
Then the following are equivalent.
\begin{itemize}
\item[1.] $\ell^1$ is not isomorphic to a subspace of $X$.
\item[2.] Every integral operator from $Y$ to $X^*$ is compact, for any $Y$.
\item[3.] Every integral operator from $\ell^1$ to $X^*$ is compact.
\item[4.] Every integral operator on $X^*$ is compact.
\item[5.] Every integral operator from $X$ to $Y$ is compact, for any $Y$.
\item[6.] Every integral operator from $X$ to $X^*$ is compact.
\item[7.] Every operator from $L^1$ to $X^*$ is Dunford-Pettis.
\item[8.] Every Dunford Pettis operator from $X$ to $Y$ is compact, for 
arbitrary $Y$.
\item[9.] Every $w^*$-compact subset of $X^*$ is the norm-closed convex 
hull of its extreme points.
\item[10.] If $K$ is a weak$^*$-compact subset of $X^*$ which isomorphically 
norms $X$, then $[K] = X^*$.
\end{itemize}
The remaining equivalences assume that $X$ is separable.
\begin{itemize}
\item[11.] Every unconditional family in $X^*$ is countable.
\item[12.] Every unconditional family in $X^*$ has cardinality less than $\G$.
\item[13.] $\L(X^*,\ell^\infty)$ has cardinality $\G$.
\item[14.] $\L(X^*)$ has cardinality $\G$.
\item[15.] $X^{**}$ has cardinality $\G$.
\item[16.] $X^{**}$ has cardinality less than $2^\G$.
\item[17.] If $K$ is a weak$^*$-compact subset of $X^*$ and $(x_n)$ is a bounded
sequence in $X$, then setting $\hat x_n(k) = k(x_n)$ for all $k$ and $n$, 
any point-wise cluster point of $(\hat x_n)$ belongs to the first Baire 
class on $K$.
\item[18.] There exists an isomorphically norming $w^*$-compact subset $K$ 
of $X^*$ so that if $(x_n)$ is as in 17, then $(\hat x_n)$ has a point-wise
cluster point which is universally measurable on $K$.
\item[19.] There exists a $K$ as in 18 so that if $(x_n)$ is as in 17, the 
cardinality of the set of point-wise cluster points of $(\hat x_n)$ on $K$, 
is less than $2^\G$.
\end{itemize}
\end{thm}


The implications $7\Rightarrow 1$, $12\Rightarrow 1$ follow quickly from a 
classical  theorem of Pe{\l}czy\'nski \cite{P}, and 
the second of these does not require the separability of $X$. 
We prove a generalization of Pe{\l}czy\'nski's result in the Appendix.

After proving Theorem~\ref{thm:main} and discussing some complements, we 
apply it in some detail to the question of when  $\ell^1$ embeds in the 
injective tensor product of two Banach spaces. 

\begin{proof} 
$1\Rightarrow 2$. 
Let $T:Y\to X^*$ be an integral operator. 
Thus there exists a probability space $(\Omega,\S,\mu)$ and 
operators $U:Y\to L^\infty (\mu)$ and $V:L^1(\mu)\to X^*$ such 
that the following commutative diagram holds
\begin{equation}\label{commdiag1}
\begin{CD}
L^\infty(\mu) @>i>> L^1(\mu)\\
@AA{U}A	@VV{V}V\\
Y @>T>> X^*
\end{CD}
\end{equation}
Here, $i:L^\infty (\mu) \to L^1 (\mu)$ denotes the canonical injection. 
Suppose $T$ is not compact. 
Then there exists a bounded sequence $(y_n)$ such that $(Ty_n)$ has no 
convergent subsequence. 
Then (as is standard), after passing to a subsequence if necessary, we 
may assume that there is a $\delta>0$ so that 
\begin{equation}\label{delta-std}
\|Ty_n - Ty_m\|  \ge \delta\ \text{ for all }\ n\ne m\ .
\end{equation}
Since $i$ is weakly compact, so is $T$; so again after passing to a 
subsequence if necessary, we may assume that $(Ty_n)$ converges weakly.
But now if we consider the sequence $(z_n)$ defined by 
\begin{equation}\label{z-seq}
z_n = y_{2n} - y_{2n-1}\ \text{ for all }\ n\ ,
\end{equation}
then $(z_n)$ is also bounded of course, and
\begin{equation}\label{z-bounded}
Tz_n\to 0\ \text{ weakly as }\ n\to\infty\ ,\ \text{ and }\ 
\|Tz_n\| \ge \delta\ \text{ for all }\ n\ .
\end{equation}
Thus it follows by the Hahn-Banach theorem that we may choose a sequence 
$(x_n)$ in the unit ball of $X$ such that 
\begin{equation}\label{delta2}
|(Tz_n) (x_n)| \ge \frac{\delta}2 \ \text{ for all }\ n\ .
\end{equation}

Now by the $\ell^1$-Theorem \cite{Ro1}, $(x_n)$ has a weak-Cauchy subsequence, 
so let us just assume that $(x_n)$ is itself weak-Cauchy (note that 
obviously \eqref{delta2} holds for subsequences
$(z'_n,x'_n)$ of the pair $(z_n,x_n)$).

But now it follows that we may choose $n_1<n_2< \cdots$ such that 
\begin{equation}\label{delta4}
|Tz_{n_{k+1}} (x_{n_k})| < \frac{\delta}4\ \text{ for all }\ k\ .
\end{equation}
Indeed, let $n_1=1$; since $Tz_n\to 0$ $w^*$, we may choose $n_2 >n_1$ 
such that $|Tz_{n_2} (x_{n_1})| <\frac{\delta}4$.
Having chosen $n_k$, choose $n_{k+1} >n_k$ such that \eqref{delta4} 
holds. 
But now it follows that 
\begin{equation}\label{eq7}
|Tz_{n_{k+1}} (x_{n_{k+1}} - x_{n_k}) | > \frac{\delta}4\ \text{ for all }\ 
k\ ,
\end{equation}
and
\begin{equation}\label{eq8}
x_{n_{k+1}} - x_{n_k} \to 0\ \text{ weakly.}
\end{equation}
Thus finally, after just re-lettering everything 
we have that 
\begin{align}	
&(z_n)\ \text{ is a bounded sequence in $Y$, $(x_n)$ and $(Tz_n)$ both}
\label{eq9}\\ 
&\text{converge to zero weakly, and (5) holds for some }\delta >0. \notag
\end{align}

Now regarding $X\subset X^{**}$, we may thus write that 
\begin{equation}\label{5dup}
|\langle Tz_n,x_n\rangle | = |\langle  T^* x_n,z_n\rangle | \ge \frac{\delta}2 
\ \text{ for all }\ n\ .
\end{equation}
But $T^*$ is also integral, and in fact admits the factorization
\begin{equation}\label{6dup}
\begin{CD}
L^\infty (\mu) @>i>> L^1 (\mu)\\
@AA{V^*}A 	@VV{\tilde U}V\\
X^{**} @>{T^*}>> Y^*
\end{CD}
\end{equation} 
where $\tilde U = U^* \mid L^1 (\mu)$, $L^1(\mu)$ regarded as contained in 
$L^1 (\mu)^{**} = (L^\infty (\mu))^*$. 

But finally, since $L^\infty (\mu)$ has the Dunford-Pettis property and 
$\tilde U i$ is weakly compact and $V^* x_n\to 0$ weakly,
\begin{equation}\label{7dup}
\|T^* x_n\|= \| \tilde u iV^* (x_n)\| \to 0\ \text{ as }\ n\to\infty\ .
\end{equation}
Of course since $(z_n)$ is bounded, this contradicts \eqref{5dup}.

$2\Rightarrow 3$, $2\Rightarrow 4$, and $2\Rightarrow 6$ are trivial.

Next, we show that 
$\text{not }1\Rightarrow \text{not }3$, $\text{not }4$, and $\text{not }6$, 
establishing the equivalence of 1, 2, 3, 4, and 6.
We shall need the following basic fact concerning integral operators:
{\em If $Y\subset Z$, $B$ is complemented in $B^{**}$,  
and $T$ is an integral operator from $Y$ to $B$, then 
$T$ extends to an integral operator from $Z$ to $B$.}
To see this, choose a probability measure space $(\Omega,\S,\mu)$ and 
operators $U: Y\to L^\infty (\mu)$ and $V: L^1 (\mu)\to B$ so that 
the following diagram commutes:
\begin{equation*}
\begin{CD}
L^\infty (\mu) @>>i> L^1 (\mu)\\
@AA{U}A		@AA{V}A\\
Y @>>T>  B			
\end{CD}
\end{equation*}
(where 		
$i$ is the canonical ``identity'' map).
Then since $L^\infty (\mu)$ has the Hahn Banach extension property, we may 
choose an operator $\tilde U: Z\to L^\infty (\mu)$ extending  $U$, (with 
$\|\tilde U\| = \|U\|$).
Thus $Vi\tilde U$ is an integral operator extending $T$, proving our assertion.

Now suppose $\ell^1$ embeds in $X$. 
Then by a theorem of Pe{\l}czy\'nski \cite{P} 
\begin{equation}\label{8dup}
(C([0,1]))^*\ \text{ embeds in }\ X^*\ .
\end{equation}
(See the Appendix for a more general result.)

Thus in particular, 
\begin{equation}\label{9dup}
L^1\ \text{ embeds in }\ X^*\ .
\end{equation}

Let $Q:\ell^1 \to C[0,1]$ be a quotient map, and let $i: C[0,1])\to L^1$ 
be the canonical ``identity'' map and also let $U:L^1\to X^*$ be an 
isomorphic embedding.

Then $T = UiQ$ is a non-compact integral operator from $\ell^1$ to $X^*$, 
proving not~3.
But if we let $Y$ be a subspace of $X$ and $V:Y\to \ell^1$ a surjective 
isomorphism, then the map $S = TV$ is a non-compact integral operator 
from $Y$ to $X^*$, and so has an integral operator extension from $X$ 
to $X^*$, proving not~6. 
Finally, let $Z$ be a subspace of $L^1$ isometric to $\ell^1$, and 
$A:Z\to\ell^1$ a surjective isometry.
Then $TA$ is a non-compact integral operator from $Z$ to $X^*$, and so has 
an integral operator extension from $X^*$ to $X^*$ by the basic fact above, 
proving not~4.

$2\Rightarrow 5$. 
Let $T: X\to Y$ be an integral operator. 
Then $T^* :Y^* \to X^*$ is also integral. 
Hence $T^*$ is compact, so $T$ is compact.
(The result that $1\Rightarrow 5$ follows from a result due to Pisier; 
see remark~4 below.)

$5\Rightarrow 6$ is trivial.

$1\Rightarrow 7$. 
Suppose to the contrary that $T:L^1\to X^*$ is a non Dunford-Pettis operator.
It follows that we may choose a sequence $(f_n)$ in $L^1$ so that 
$f_n\to 0$ weakly but $\|Tf_n\| \not\to 0$.
Therefore we may choose a subsequence $(f'_n)$ of $(f_n)$ so that for 
some $\delta >0$, 
\begin{equation}\label{nonDP}
\|Tf'_n\| > \delta \ \text{ for all }\ n\ .
\end{equation}
For each $n$, choose $x_n\in X$ with $\|x_n\| =1$ and 
\begin{equation}\label{nonDP-n}
|\langle Tf'_n ,x_n\rangle | > \delta\ .
\end{equation}
By the $\ell^1$-Theorem, by passing to a  further subsequence if 
necessary, we may assume that $(x_n)$ is a weak-Cauchy sequence.
But then it follows that 
\begin{equation}\label{weak-Cauchy}
(T^* x_n) \ \text{ is a weak Cauchy sequence in }\ L^\infty
\end{equation} 
(where we regard $X$ as canonically embedded in $X^*$).

Since $L^1$ has the Dunford-Pettis property, it follows that 
\begin{equation}\label{DP-property}
|\langle f'_n, T^* x_n\rangle |\to 0\ \text{ as }\ n\to\infty\ ,
\end{equation}
which contradicts \eqref{nonDP-n}. 
(A Banach space $B$ has the Dunford-Pettis property if every weakly compact 
operator $T:B\to Y$ is Dunford-Pettis, for all $Y$. 
By fundamental results of Grothendieck \cite{G1}, this is equivalent to: 
$(b_n)$ weakly null in $B$,
$(f_n)$ weakly null in $B^*$ implies $f_n(b_n)\to 0$ as $n\to\infty$. 
It is then a standard exercise to show that in fact if $(b_n)$ is weakly 
null and $(f_n)$ is weak-Cauchy, still $f_n(b_n)\to 0$ as $n\to\infty$.)

Not 1$\Rightarrow $ not 7.
Since $L^1$ is isomorphic to a subspace of $X^*$ when 1 fails, this is 
immediate: 
An (into) isomorphism $T:L^1\to X^*$ is obviously not Dunford-Pettis.

$1\Rightarrow 8$.
Let $(x_n)$ be a bounded sequence in $X$. 
Choose (by the $\ell^1$-Theorem) $(x'_n)$ a weak Cauchy subsequence of $(x_n)$;
then given $(n_i)$ and $(m_i)$ strictly increasing
sequences of positive integers, 
$(x'_{n_i} - x'_{m_i})$ is weakly null,  and hence by hypothesis
$$\|T(x'_{n_i} - x'_{m_i})\| = \|Tx'_{n_i} - T'x _{m_i}\| \to 0
\ \text{ as }\  i\to\infty\ .$$
This implies $(Tx'_n)$ is a Cauchy sequence in $Y$, so it converges since 
$Y$ is complete; thus $T$ is compact.

$8\Rightarrow 5$.
Integral operators are Dunford-Pettis operators because $L^\infty (\mu)$-spaces
have the Dunford-Pettis property and integral operators factor through 
the ``identity'' map $i:L^\infty (\mu) \to L^1 (\mu)$, for some probability 
measure $\mu$; of course $i$ is weakly compact.

$1\Rightarrow 9$. 
Suppose to the contrary that $(u_\alpha)_{\alpha <w_1}$ is an uncountable 
unconditional family in $X^*$; let $U$ be the norm closure of its linear span.
Assume $\|u_\alpha\| = 1$ for all $\alpha$; letting $(u_\alpha^*)$ be the 
functions in $U^*$ biorthogonal to $(u_\alpha)$, choose $K$ so that 
$\|u_\alpha^K\| \le K$ for all $\alpha$.
Now we have the following fundamental claim:
\begin{equation}\label{eq:fundamental}
\text{Given $f\in U^*$, then $W_f = \{\alpha :f(u_\alpha) \ne 0\}$ is 
countable.}
\end{equation}
If this were false, say then $f= U^*$, $\|f\| =1$, and $W_f$ is countable. 
Then we may pass to an uncountable subset $\Gamma$ of $W_f$ such that there 
exists a $\delta >0$ so that 
\begin{equation}\label{fundamental-false}
|f(u_\alpha)| \ge \delta\ \text{ for all }\ \alpha\in\Gamma\ .
\end{equation} 
But now a standard argument shows that 
\begin{equation}\label{std-argument}
(u_\alpha)_{\alpha\in\Gamma} \ \text{ is equivalent to the natural basis of }\ 
\ell^1 (\Gamma)\ .
\end{equation} 
Thus $\ell^1 (\Gamma)$ embeds in $X^*$, which implies $\ell^1$ embeds in $X$,
by a result of Pe{\l}czy\'nski \cite{P}, Hagler \cite{H1},  
a contradiction.

$1\Rightarrow 9$. 
This is due to R. Haydon \cite{Hay}. 
For $X$ separable, this had previously been proved by E. Odell and myself 
\cite{OR}.

$9\Rightarrow 10$. 
If $K$ satisfies the hypothesis of 10, so does $\tilde K = \{\alpha k : 
\alpha$ is a scalar, $|\alpha| =1$, $k\in K\}$. 
But then it follows that the $w^*$-closed convex hull $W$ of $\tilde K$ has 
non-empty interior, and hence since then also $\text{Ext }W\subset \tilde K$ 
and 9 implies $W$ is the norm closed convex hull of $\text{Ext }W$, 
$[K] = [\tilde K] =X^*$.

$10\Rightarrow 1$. 
This follows from a result of G.~Godefroy \cite{G1}.
Indeed, assume that $\ell^1$ embeds in $X$.
Then it is proved in \cite{Gr1} that there exists an equivalent norm 
$\trinorm{\cdot}$ on $X$ such that if $K$ denotes the $w^*$ closure of 
the extreme points of the ball of $X^*$ in the dual norm induced by 
$\trinorm{\cdot}$, then $[K] \ne X^*$.
But of course then $K$ is an isomorphically norming $w^*$-compact 
subset of $X^*$ in its original norm, proving that 10 does not hold.
There is a minor point here that requires explanation, however. 
The cited result  of Godefroy's requires the fundamental case where 
$X=\ell^1$ itself, for the proof.
But the argument given in \cite{G1} is only valid for the case of complex 
scalars.
The result for real scalars may be deduced from the work  in \cite{G1} 
as follows: 
First, to avoid ambiguity, let $\ell_\real^p$, resp. $\ell_\complex^p$, 
denote $\ell^p$ for real scalars, resp. for complex scalars, $p=1$ 
or $\infty$.  
It is proved in \cite{G1} that for complex scalars, there exists an 
equivalent norm $\trinorm{\cdot}$ on $\ell_\complex^1$ so that if $K$ is 
as above, and $Y$ is the closed linear span of $K$ over the complex 
scalars, then there exists an infinite subset $M$ of $\nat$ with infinite 
complement so that 
\begin{equation}\label{eq:10to1}
\inf \{ |y(n) - y(n')| : n\in M, \ n'\in N \sim M\} =0 \ \text{ for all }\ 
y\in Y\ .
\end{equation}

Now let $\trinorm{\cdot}^*$ be the dual norm induced on $\ell_\complex^\infty$,
and just regard $(\ell_\complex^\infty,\trinorm{\cdot}^*)$ as a real 
Banach space. 
Now if we take the standard norm on $\ell_\complex^\infty$ and regard this 
as a real Banach space, we obtain $\ell_\real^\infty \oplus \ell_\real^\infty$
under the norm 
\begin{equation}\label{stdnorm}
\| (a_j) \oplus (b_j)\|= \sup_j \sqrt{a_j^2 + b_j^2}\ \text{ for }\ 
(a_j),(b_j) \in \ell_\real^\infty
\end{equation}
which is obviously equivalent to the standard norm on $\ell_\real^\infty \oplus 
\ell_\real^\infty$, which of course is isometric to $\ell_\real^\infty$. 
Then it follows also that the norm $\trinorm{\cdot}^*$ must be equivalent 
to the standard norm $\|\cdot\|_\infty$ on $\ell_\real^\infty \oplus 
\ell_\real^\infty$ (where obviously we take the isomorphism 
$(a_j +ib_j) \to ((a_j),(b_j))$ for $(a_j +ib_j) \in \ell_\complex^\infty$).

Now let $c_{00}$ denote the space of all sequences of reals which are 
ultimately zero.
Define a norm $\trinorm{\cdot}$ on $c_{00} \oplus c_{00}$ by 
\begin{equation}\label{ultimatezero}
(a_j) \oplus (b_j) = \sup \Big\{ |\sum (\alpha_j a_j + \beta_j b_j)| 
: \trinorm{(\alpha_j) \oplus (b_j)}^* =1\Big\}\  .
\end{equation}
It follows easily that $\trinorm{\cdot}$ is equivalent to the $\ell^1$-norm 
on $\ell_\real^1 \oplus \ell_\real^1$. 
Moreover, we have that a bounded sequence $(f_n)$ in $\ell_\complex^\infty$ 
converges in the $w^*$ topology on $\ell_\complex^\infty$ induced by 
$\ell_\complex^1$ iff  
\begin{equation}\label{f-sub-n1}
\lim_{bn\to\infty} f_n (j)\ \text{ exists for all }\ j
\end{equation}
iff 
\begin{equation}\label{f-sub-n2}
\lim_{n\to\infty} \Re f_n(j)\ \text{ and }\ \lim_{n\to\infty}\Im f_n(j) 
\ \text{ exist for all }\ j
\end{equation}
iff 
\begin{equation}\label{f-sub-n3}
\{ (\Re f_n \oplus \Im f_n)_{n=1}^\infty\ \text{ converges in the 
$w^*$-topology on }\ \ell_\real^\infty \oplus \ell_\real^\infty\}\ .
\end{equation}
It follows that the $w^*$-topology on $\ell_\complex^\infty$ is the same
as that on $\ell_\real^\infty \oplus \ell_\real^\infty$, and hence the 
unit ball of $(\ell_\real^\infty \oplus \ell_\real^\infty, \trinorm{\cdot}^*)$ 
is compact in the $w^*$-topology; thus by the bipolar  theorem, that 
ball is precisely the dual ball of $(\ell_\real^1 \oplus \ell_\real^1,
\trinorm{\cdot})$; hence the set $K$ defined above is exactly the same 
as that defined for the real scalars case  above. 
Finally, it follows that if $Y_\real$ denotes the norm closed linear 
span of $K$ over real scalars, then we have that if 
$(a_j)_{j=1}^\infty \oplus (b_j)_{j=1}^\infty \in Y_\real$, 
$(a_j + ib_j)_{j=1}^\infty \in Y$, and hence by \eqref{eq:10to1}, we 
obtain that 
\begin{equation}\label{normclosed}
\inf\{ \max \{ |a_n - a'_n|, |b_n-b'_n| \} : n\in M,\ n'\in \nat\sim M\} =0
\end{equation}
which obviously implies that $Y_\real \ne \ell_\real^\infty \oplus 
\ell_\real^\infty$.
Thus it follows $X= (\ell_\real^1 \oplus \ell_\real^1,\trinorm{\cdot})$ 
is isomorphic to $\ell_\real^1$ but for $K$ as defined above, the 
$w^*$-closed convex ball of $K$ equals $\Ba X^*$ but $[K]\ne X^*$, 
completing the proof for real scalars.

$11\Rightarrow 12$. 
This is trivial.

It follows from the main result in \cite{OR} that every element of 
$X^{**}$ is the limit of a weak Cauchy sequence in $X$, which yields 15. 
($1\Rightarrow 15$ is also given in \cite{OR}.)

$1\Rightarrow 13$. 
It is easily seen that 
\begin{equation}\label{card1}
\card \L (X^*,\ell^\infty) \ge \G\ .
\end{equation}
Indeed, just fix $z\in \ell^\infty$, $z\ne0$, $x\in X$, $x\ne 0$, and note that 
the operator $T$ on $X^*$ defined by $T(x^*) = x^* (x) z$ is 
thus non-zero; hence 
\begin{equation}\label{eq-new23}
r\to rT \ \text{ is a 1--1 map of }\ \real \ \text{ into  } \ 
\L (X^*,\ell^\infty)\ . 
\end{equation}
Thus it remains to prove
\begin{equation}\label{card2}
\card \L(X^*,x) \le \G\ .
\end{equation}
%
Now since $\Ba \ell^1$ is $w^*$-dense in $\Ba (\ell^\infty)^*$ by Goldstein's 
Theorem. 
then 
\begin{equation}\label{Zisomorphic}
\text{There exists a countable subset $D$ of $\Ba \ell^\infty$ 
which is weak*-dense in it.}
\end{equation}
But then it follows that if $T\in \L^1 (X^*,\ell^\infty)$, 
$T^*$ is $w^*$-continuous and is thus determined by its values on $D$. 
Thus we deduce that 
\begin{equation}\label{card4}
\card \L (X^*,\ell^\infty) \le \card (X^{**})^D = 
\G^{\aleph_0} = \G\ .  
\end{equation}

$13\Rightarrow 14$ 
is obvious since the dual of any separable Banach space
is isometric  to a subspace of $\ell^\infty$.

$14\Rightarrow 15$ is trivial 
since $X^{**}$ is isometric to a subspace of $\L(X^*)$.

$15\Rightarrow 16$ is trivial.

$1\Rightarrow 17$.
It follows by the $\ell^1$-theorem that given $(\hat x_{n_i})$ a subsequence
of $(\hat x_n)$, then $(x_{n_i})$ has a weak-Cauchy subsequence, which 
implies that $(\hat x_{n_i})$ has a subsequence pointwise convergent on $K$.
Now the results of \cite{Ro2} prove the conclusion of 17 (and also imply, by
the way, that any pointwise cluster point of $(\hat x_n)$ is the limit of 
a pointwise convergent subsequence of $(x_n)$).

$17\Rightarrow 18$.
Let $K$ be the unit ball of $X^*$. 
Then of course $(\hat x_n)$ has a pointwise cluster point on $K$, which is 
just an element of $X^{**}$ restricted to $K$. 
This is a Baire-one functions on $K$ by \cite{OR}, and so of course is 
Borel measurable and hence universally measurable.

$17\Rightarrow 19$.
The cardinality of the class of Baire-one functions on $K$ equals 
$\G < 2^\G$.

$18\Rightarrow 1$.
It suffices to prove that any bounded sequence $(x_n)$ in $X$ has a 
subsequence $(x_{n_i})$ so that $(\hat x_{n_i})$ converges pointwise on $K$. 
For then,		 
it follows by the 
Hahn-Banach, Riesz-representation, and bounded convergence theorems that 
$(x_{n_i})$ is a weak-Cauchy sequence, and so 1 holds. 
But if this is not the case, then we find a bounded sequence $(x_n)$ in $X$ 
such that $(\hat x_n)$ has no pointwise convergent subsequence on $K$. 
It now follows by a theorem of Bourgain-Fremlin-Talagrand \cite{BFT} 
that $(\hat x_n)$ has a subsequence, none of whose pointwise cluster 
points are universally measurable on $K$, thus contradicting 18.
(For an alternate proof of the cited result, see Theorem~3.18 in \cite{Ro3}.)

It remains to prove that 12, 16, and 19 imply 1.
We shall prove the contrapositive implications instead.
So we assume for the rest of the proof that $\ell^1$ embeds in $X$ 
and $X$ is separable.
(We don't need the separability assumption for the first two implications.)

$\ell_\G^1$ is isometric to the space of atomic Borel measures on $[0,1]$, 
and so isometric to a subspace of $C[0,1)^*$. 
Thus by \eqref{eq8}, 
\begin{equation}\label{atomic}
\ell_\G^1 \ \text{ is isomorphic to a subspace of $X^*$,} 
\end{equation}
so 12 does not hold.

It also follows that 16 does not hold, for by 
\eqref{atomic}, 
$\ell_\G^\infty$ is isomorphic to a quotient space of $X^{**}$, and hence
\begin{equation}\label{ell-isomorphic}
\card X^{**} \ge \G^\G = 2^\G\ .
\end{equation}
(These implications are of course known.)

Not $1\Rightarrow$ Not 19.
Let $K$ be a $w^*$-compact isomorphically norming subset of $X^*$. 
Then it follows again by the theorems cited in the proof of $18\Rightarrow 1$,
that if $(x_n)$ is a bounded sequence in $X$ such that  $(\hat x_n)$ 
converges pointwise on $K$, then $(x_n)$ is a weak Cauchy sequence in $X$.
Thus there must exist some bounded sequence $(x_n)$ in $X$ such that 
$(\hat x_n)$ has no pointwise convergent subsequence on $K$. 
Hence this implication follows (since obviously the family of pointwise 
cluster  points on $K$ has cardinality at most $2^\G$) from 
\renewcommand{\qed}{}
\end{proof}

\begin{lem}\label{lem2}
Let $K$ be a compact Hausdorff space and $(f_n)$ be a bounded sequence of 
continuous scalar-valued functions on $K$ which has no pointwise convergent 
subsequence. 
Then the family $F$ of pointwise cluster points of $(f_n)$ on $K$ has 
cardinality at least $2^\G$.
\end{lem}

\begin{proof}
We may obviously assume the $f_n$'s are real valued since either the real or
imaginary parts of the $f_n$'s have no pointwise convergent subsequence, and 
the real or imaginary part of a pointwise cluster point of $(f_n)$ is also 
a pointwise cluster point of the real or imaginary parts of the $f_n$'s.
By the proof of the $\ell^1$-Theorem (\cite{Ro1}; see also \cite{Ro3}), 
there exists a subsequence $(f'_n)$ of $(f_n)$ with the following property:
\begin{quote}
$(*)$ There exist real numbers $r$ and $\delta$ with $\delta>0$ such that 
setting $A_n = \{k: f'_n(k) \le r-\delta\}$ and 
$B_n = \{k :f'_n(k) \ge r+\delta\}$, then these  are non-empty sets for 
all $n$, such that  $((A_n,B_n))$ is a Boolean independent sequence; 
that is, if one sets $+A_n = A_n$ and $-A_n =B_n$, then for any infinite 
sequence $(\E_j)$ with $\E_j = \pm1$ for all $j$,
\end{quote}
\begin{equation}\label{Boolean1}
\bigcap_{j=1}^n \E_j A_j\ \text{ is non-empty for all }\ n\ .
\end{equation}
This implies
\begin{equation}\label{Boolean2}
\bigcap_{j=1}^\infty \E_j A_j\ \text{  is non-empty,}
\end{equation}
since for all $n$, $f'_n$ is continuous, and thus $A_n,B_n$ are closed 
and non-empty, and so \eqref{Boolean2} follows by the compactness of $K$.
Now let $\U$ be the family of all non-principal ultrafilters on $\nat$ 
(cf. \cite{CN} for the definition and standard properties of 
ultrafilters). 
Then as is classical, $U$ may be identified the $\beta \nat\sim \nat$, where 
$\beta \nat$ denotes the Stone-C\v each compactification of $\nat$, and so by 
a classical theorem in topology (cf. Theorem~2, page 132 of \cite{E}), 
\begin{equation}\label{betaW}
\card \U = 2^\G\ .
\end{equation}
(This is also explicitly given in 7.4 Corollary, page 146 of \cite{CN}.)

For each $U\in \U$, define a function $f_U$ on $K$ by 
\begin{equation}\label{u-in-U}
f_U (k) = \lim_{n\in U} f'_n(k)\ \text{ for all }\ k\in K\ .
\end{equation}
Then as is standard, $f_U$ is a pointwise cluster point of $(f'_n)$ and 
hence of $(f_n)$. 
Thus to complete the proof of the Lemma, it suffices to show that 
\begin{equation}\label{lemma-proof}
f_U \ne f_V\ \text{ if }\ u\ne v\ ,\quad U,V\in \U\ .
\end{equation}
Given $U\ne V$ in $u$, choose $M$ an infinite subset of $\nat$ 
such that $M \in U$, $M \notin V$.  
Then as cofinite sets belong to any non-principal ultrafilter, 
\begin{equation}\label{cofinite-sets}
\nat \sim M \defeq L\ \text{ is infinite and thus belongs to }\ V\ .
\end{equation}
(The latter statement holds since $M$ is infinite and $M\notin V$.) 
Now define $(\E_j)$ by 
\begin{equation}\label{defineE}
\E_j = 1\ \text{ if }\ j\in M\ \text{ and }\ \E_j =-1\ \text{ if} \ j\in L\ .
\end{equation}
Let $k$ be a point in $\bigcap_{j=1}^\infty \E_j A_j$ (such exists by 
\eqref{Boolean2}). 
Thus by definition of $((A_n,B_n))$ 
\begin{equation}\label{Asubn-Bsubn}
f'_n(k) \le r-\delta \ \text{ for all }\ n\in M\ .
\end{equation}
Now it follows that if $U_M = \{M \cap A: A\in U\}$, then 
$U_M$ is a non-principal ultrafilter on $M$, and 
\begin{equation} \label{eq45}
f_U (k) = \lim_{k\in u_\mu} f'_n(k)\ .
\end{equation}
Thus it follows by \eqref{Asubn-Bsubn} that 
\begin{equation}\label{eq46}
f_U (k) \le r-\delta\ .
\end{equation}
By exactly the same reasoning, we obtain that 
\begin{equation}\label{eq47}
f_V (k) \ge r+\delta\ .
\end{equation}
Thus \eqref{eq46} and \eqref{eq47} show $f_U\ne f_V$, completing the
proof of the Lemma, and thus the proof of Theorem~\ref{thm:main}.
\end{proof}

\subsection*{Remarks} $\quad$

%

1a. The equivalence we have used in the proof of $10\Rightarrow 1$ is 
actually quantitative. 
That is, we have the following fact.

\begin{propA}\label{quantitative}
Given $X$ and $\lambda \ge 1$, the following are equivalent:
\begin{itemize}
\item[(i)]  There exists a weak*-compact $\lambda$-norming subset $K$ 
of $\Ba X^*$ so that $[K]\ne X^*$.
\item[(ii)] There exists a norm $\trinorm{\cdot}$ on $X$ so that 
$$\| \cdot\| \le \trinorm{\cdot} \le \lambda \|\cdot\|\qquad 
(\|\cdot\|\ \text{ the original norm})$$
with $[K]\ne X^*$, where $K$ is the $w^*$-closure of the extreme points 
of $\Ba (X^*,\trinorm{\cdot}^*)$.
\end{itemize}
\end{propA}

\begin{proof}
(i) $\Rightarrow$ (ii): 
Let $W = \{\alpha k : |\alpha|=1$, $k\in K\}$. 
Then $W$ is also $w^*$-compact. 
Then it follows by the geometrical form of the Hahn-Banach Theorem that 
if $\tilde K$ denotes the $w^*$-closed convex hull of $W$, then 
$$\tilde K \subset \Ba X \subset \lambda \Ba \tilde K\ ;$$
in turn, we then easily obtain a norm $\trinorm{\cdot}$ on $X$ 
satisfying the inequality in (ii) such that 
$\tilde K = \Ba (X^*,\trinorm{\cdot}^*)$. 
But then it follows that $\Ext \tilde K\subset W$ which implies (ii).
(ii) $\Rightarrow$ (i) is immediate for $K$ is then $\lambda$-norming.
\end{proof}

Now the arguments in \cite{G1} yield 

\begin{propB}\label{propB}
There exists an absolute constant $C$ so that if $\ell^1$ embeds in $X$, 
there exists a $C$-norming $w^*$-compact subset $K$ of $\Ba X^*$ with 
$[K]\ne X^*$.
\end{propB}

The question then arises:

\begin{prob}
If $\ell^1$ embeds in $X$, is it so that given $\ep>0$, there exists 
a $\lambda \le 1+\ep$ satisfying (i) of Proposition A?
\end{prob}

The delicate nature of the proof of 10 $\Rightarrow 1$ leads me to 
conjecture that the answer is negative for $X= \ell^1$ itself. 
Of course if so, it is natural to ask:
What is the optimal value of $C$ in Proposition~B?
Now of course if $X = C(K)$ for some uncountable compact metric space, 
$\lambda =1$  works.
However $\lambda =1$ does not work for $X=\ell^1$. 
For then by the arguments sketched above, we  would have that the norm-closed 
linear span of the $w^*$-closure of $\Ba \ell^\infty$ would be unequal to 
$\ell^\infty$. 
But standard arguments show that the norm-closed convex hull of 
$\Ext \Ba \ell^\infty = \Ba \ell^\infty$. 

1b.  Some interesting equivalences  are also obtained in \cite{G1}, 
complementary to the above discussion. 
The following notion is introduced there:

\begin{Define} 
Let $K$ be a non-empty closed bounded convex subset of $X^*$, and 
$W\subset K$. 
$W$ is called a {\em boundary}
of $K$ if for all $x\in X$, $\sup [\Re (\chi x)\mid W]$
is attained on $W$.
\end{Define}

The following result  is obtained in \cite{G1}, generalizing the 
theorem of \cite{OR} that $9\Rightarrow 1$ for separable $X$, and 
yielding an analogy of James' famous characterization of weakly 
compact convex sets.

\begin{Thm}
Assume $X$ is separable. 
Then the following are equivalent.
\begin{itemize}
\item[(i)]  $\ell^1$ does not embed in $X$.
\item[(ii)] for all $w^*$-compact subsets $K$ of $X^*$, if $W$ is a 
boundary of $K$, then $K$ is the norm-closed convex hull of $W$.
\item[(iii)] if $K$ is a closed bounded convex subset of $X^*$ which 
is a boundary for itself, then $K$ is $w^*$-compact.
\end{itemize}
\end{Thm}

2. I am indebted to Welfeng Chen for a stimulating conversation 
concerning the implication $1\Rightarrow 12$.

3.  Suppose $X$ is separable and $X^*$ is non-separable. 
A remarkable result of Stegall yields that then $X^*$ has a biorthogonal 
family of cardinality the continuum \cite{S}. 
Since there are now known many separable spaces $X$ not containing
$\ell^1$ with $X^*$ nonseparable 
(cf.\ \cite{AMP}, \cite{H2}, \cite{J3}, \cite{K}, \cite{LS}, \cite{Ro5}), 
we cannot significantly weaken the unconditionality assertion in 11 and 12. 
This same result incidentally shows that $\dens X^{**} = \dens X^* =\G$;
thus assuming $\ell^1$ does not embed in $X$, we obtain that all 
cardinal measures of the size of $X^*$, $X^{**}$, $\L(X^*)$ yield the 
same result.
(For a metric space $M$, $\dens M$ denotes the least cardinality of 
a dense subset.)

%

4. The following essentially known result gives equivalences for the 
embedability of $\ell^1$ in $X$ and the structure of $p$-absolutely summing
operators on $X$, analogous to the equivalences 1.\ through 4.\ of 
Theorem~\ref{thm:main}.

\begin{Thm}\label{thm:summing}
Let $X$ be given.
The following are equivalent.
\begin{itemize}
\item[1.] $\ell^1$ does not embed in $X$.
\item[2.] For all $Y$ and $1\le p<\infty$, every $p$-absolutely summing 
operator from $X$ to $Y$ is compact.
\item[3.] Every $2$-summing operator from $X$ to $\ell^2$ is compact.
\item[4.] Every $2$-summing operator from $X$ to $X^*$ is compact.
\end{itemize}
\end{Thm}

\begin{proof} 
$1\Rightarrow 2$.
This is due to G.~Pisier (Corollary 1.7, part (iii) of \cite{Pi}). 
For the sake of completeness, we give the argument.
Let $T:X\to Y$ be a $p$-absolutely summing operator.
By a fundamental theorem of 
Pietsch (cf., \cite{LP} or  Theorem~1.3 of \cite{Pi}), 
there  exists a regular Borel probability measure $\mu$ on $K=\Ba X^*$ 
endowed with the $w^*$-topology and a $C<\infty$ so that 

(a) $\|Tx\| \le C(\int_K |w(x)|^p\ell\mu (w))^{1/p}$ for all $x\in X$.

\noindent
Now let $(x_n)$ be a bounded sequence in $X$, and choose a weak-Cauchy
subsequence $(x'_n)$ by the $\ell^1$-Theorem. 
We claim that 

(b) $(Tx'_n)$ converges in the norm topology of $Y$.

\noindent
To prove this, since $Y$ is a Banach space, we employ the following 
elegant characterization of Cauchy sequences, due to Pe{\l}czy\'nski:
It suffices to show that 

(c) $\|Tx'_{n_i} - Tx'_{m_i}\| \to 0$ as $i\to\infty$ for all 
strictly increasing sequences $(n_i),(m_i)$ of $\nat$.

\noindent
But given such sequences, $(x_{n_i}-x_{m_i})$ is a weakly null sequence in $X$. 
Now letting $i:C(K)\to L^p(\mu)$ be the natural injection and $U:X\to C(K)$
the canonical map given by: $(Ux)(k)=k(x)$ for all $x\in X$, $k\in K$, then 
since $i$ is weakly compact and $C(K)$ has the Dunford-Pettis property,
\begin{itemize}
\item[(d)] $\|iU(x'_{n_i} - x'_{m_i})\| \to 0 \ \text{ as }\ i\to\infty$
\end{itemize}
But (a) yields that 
\begin{itemize}
\item[(f)] $\|T (x'_{n_i} - x'_{m_i})\| \le C\| iU (x'_{n_i} - x'_{m_i}\|
_{L^p(\mu)} \ \text{ for all }\ i$.
\end{itemize}
\renewcommand{\qed}{}
\end{proof}

$2\Rightarrow 3$, $2\Rightarrow 4$
are trivial. 
Now suppose $\ell^1$ embeds in $X$.
Thus by Pe{\l}czy\'nski's theorem, \eqref{eq8} holds, and since $\ell^2$ 
embeds in $L^1$, which thus embeds in $X^*$, it is obviously enough to show 
that condition~3 of the Theorem fails to hold. 
Now simply let $T:\ell^1\to \ell^2$ be defined by:
$Te_j = b_j$ for all $j$, where $e_j$ is the $\ell^1$-basis, 
$(b_j)$ the $\ell^2$-basis. 
It follows that $T$ is absolutely summing, by Grothendieck's fundamental 
theorem \cite{G2} (although the direct proof of this elementary fact is 
much simpler).
Let $Z$ be a subspace of $X$ isomorphic to $\ell^1$ and let $S:Z\to\ell^1$ 
be a surjective isomorphism. 
It follows that $TS$ is absolutely summing. 
Hence $TS$ is 2-absolutely summing, which implies $TS$ is 2-integral and 
hence $TS$ entends to a 2-integral operator $V:X\to \ell^2$. 
$V$ is thus 2-absolutely summing, but since $T$ is not compact, 
neither is $V$.\qed
\medskip

\noindent {\it Comment\/}.
It thus follows that $1\Rightarrow 5$ of Theorem~1 can be deduced from 
the above Theorem, since integral operators are asbolutely summing.

5. The proof of $6\Rightarrow 1$ yields an integral non-compact operator 
from $\ell^1$ to $L^1$ (and also $6\Rightarrow 1$ follows from this and 
Pe{\l}czy\'nski's theorem cited there). 
The following is a more natural  example of such an operator. 
Letting $(e_n)$ be the $\ell^1$ basis, define $T$ by:
$Te_n = \sin 2\pi nx$ for all $n$.
$T$ is obviously not compact. 
To see that $T$ is an integral operator, define $S:C[0,1]\to c_0$ by:
$(Sf)_n = \int_0^1 f(x) \sin 2\pi nx\,dx$ for all $n\in N$.
Then $S$ is an integral operator, because if we define $V:L^1\to c_0$ by : 
$(Vf)_n = \in_0^1 f(x)\sin 2\pi nx$ for all $n\in N$, then $S= Vi$, where 
$i:C[0,1]\to L^1$  is the canonical map. 
Therefore $S^*$ is integral; it is easily verified that $S^* =T$.


If $K$ satisfies the hypothesis of 9, then if we renorm $X$ by 
$\trinorm{x} = \sup \{ |k(x)| :k\in K\}$, $K$ now isometrically norms
$(X,\trinorm{\cdot})$, and so also $\Omega = \{\alpha k:\alpha$ is a scalar, 
$|\alpha| =1$, $k\in K\}$ isometrically norms $(X,\trinorm{\cdot})$ and is 
$w^*$-compact.
Hence the $w^*$-closed convex hull of $\Omega$ equals $D$, the unit ball 
of $X^*$ in $\trinorm{\cdot}^*$, and moreover the extreme points of $D$ 
are contained in $\Omega$. 
Thus by 9, the norm-closed convex hull of 
$\Omega = D$, so of course $[\Omega] = [K] = X^*$. 


7. Results of Bourgain in \cite{Bo1}, \cite{Bo2} yield  the following 
remarkable result.

\begin{Thm}\label{thm:bourgain}
Assume $X$ is separable. 
Then the following are equivalent.
\begin{itemize}
\item[1.] $\ell^1$ does not embed in $X$.
\item[2.] Any closed bounded convex subset of $X^*$ with the Radon-Nikodym
property (RNP) is separable.
\item[3.] Any subspace of $X^*$ with the RNP is separable.
\item[4.] There does  not exist a compact subset $K$of $E\Ba X^*$, the 
extreme points of $\Ba X^*$, such that $K$ is homeomorphic to the Cantor set
and the canonical map $X\to C(K)$ is surjective.
\item[5.] There does not exist a subset of the extreme points of 
$\Ba X^*$ equivalent to the basis
of $\ell_\G^1$.
\end{itemize}
\end{Thm}

\noindent
(Note that $3\Rightarrow 1$ follows from Pe{\l}czy\'nski's theorem and does 
not require the separability of $X$, for $\ell_\G^1$ has the RNP.) 
He then deduces in \cite{Bo2} that if $X$ is non-separable, $E\Ba X^*$ 
is non-separable in norm. 
($EW$ denotes the extreme points of $W$.)
For if $\ell^1$ is not isomorphic to a subspace of $X$, $X^*=[E\Ba X^*]$ by
the result of Haydon \cite{Hay}.
If $\ell^1$ embeds in $X$, then there is a separable subspace $Y$ of $X$ 
so that $\card E\Ba Y^* = \G$ by $4\Rightarrow 1$ of the above Theorem, but 
as is standard, every extreme pont of $\Ba Y^*$ lifts to an extreme point 
of $\Ba X^*$, so $\card E\Ba X^* \ge \G$. 
The Theorem established by Bourgain in \cite{Bo1} which gives 
$4\To 1$, also yields the following striking improvement of Pe{\l}czy\'nski's 
theorem cited in the proof of $12\To1$ of our Theorem~\ref{thm:main}. 
{\em If $\ell^1$ embeds in $X$, $X$ separable, then for all $1>\ep>0$, there 
exists a subset $K$ of $E\Ba X^*$ homeomorphic to the Cantor set, such that 
for all $f\in (1-\ep)\Ba C(K)$, there exists an $x\in X$ with $\|x\| <1$ and 
$f(k)=k(x)$ for all $k\in K$.} 
It follows that 
{\em for all $\ep>0$, $E\Ba X^*$ has a subset $(1+\ep)$-equivalent 
to the basis for $\ell_\G^1$.}

8. Theorem 3.18 of \cite{Ro3} yields the following result, strengthening 
the result in \cite{BFT} 
used in proving $18\To 1$ of Theorem~\ref{thm:main}.

\begin{Thm}\label{thm:BFT}
Let $K$ be a compact metric space and $(f_n)$ a bounded sequence in $C(K)$ 
with no pointwise convergent subsequence. 
Then there exists a subsequence $(f'_n)$ of $(f_n)$ and a Borel probability 
measure $\mu$ on $K$ such that no point-wise cluster point of $(f'_n)$ is 
$\mu$-measurable.
\end{Thm}

It follows that we may replace the implications $18, 19\To 1$ by the 
following stronger statement:
{\em If $\ell^1$ embeds in $X$ (assumed separable) and $K$ is a weak*-compact
norming subset of $X^*$, there exists a bounded sequence $(x_n)$ in $X$ 
such that $(\hat x_n)$ {\rm (as defined in 17)} has $2^\G$ pointwise
cluster points on $K$, none of which are universally measurable.}

Applying Lemma~\ref{lem2}, we obtain from the Theorem

\begin{corA} 
If $(f_n)$ satisfies the hypotheses of the Theorem, there exists a Borel
probability measure $\mu$ on $K$ such that there is a set $W$ of cardinality 
$2^\G$, consisting of point-wise cluster points of $(f_n)$, none of which 
are $\mu$-measurable.
\end{corA}

Here is another applcation, due to Stegall.

\begin{corB}
Let $W$ be a weakly pre-compact subset of a Banach space. 
Then the closed convex hull of $W$ is weakly pre-compact.
\end{corB}

\noindent 
($W$ is called weakly pre-compact of every sequence in $W$ has a weak-Cauchy
subsequence.)

\begin{proof}
It obviously  suffices to prove that $\co W$, the convex hull of $W$, is 
weakly  pre-compact.
Any sequence in $\co W$ is contained in the convex hull of a countable subset
of $W$, so we may assume w.l.g.\ $W$ is separable.
Note incidentally  that any weakly pre-compact set is bounded, since 
weak-Cauchy sequences are bounded by the uniform boundedness principle.
Thus, suppose that  $C<\infty$, $\|w\|\le C$ for all $w\in W$, and to the 
contrary, there exists a sequence $(f_n)$ in $\co W$ with no weak-Cauchy 
subsequence.
Now letting $K = (\Ba X^*, w^*)$ and $\hat x(k) = k(x)$ for all $x\in W$ and 
$k\in K$, $(\hat f_n)\subset C(K)$ thus has no pointwise convergent subsequence 
on $K$, and so by the above Theorem, choose $\mu$ a Borel probability measure 
on $K$ and $(f'_n)$ a subsequence of $(f_n)$ such that no point-wise cluster 
point of $(\hat f'_n)$ is $\mu$-measurable.
Now letting $i:C(K)\to L^1 (\mu)$ be the canonical map, it follows by the 
bounded convergence theorem that $i(W)$ is a relatively compact subset of 
$L^1(\mu)$.
But then $coi(W) = icoW$ is also a relatively compact set. 
So then, choose $(f''_n)$ a subsequence of $(f'_n)$ such that 
$(i\widehat{f''_n})$ converges in $L^1 (\mu)$, and finally choose 
$(\tilde f_n)$ a subsequence of $(f''_n)$ such that $(f_n)$ converges 
$\mu$-almost everywhere, to a function $g$ in $L^1(\mu)$. 
But then any point-wise cluseter point of $(\widehat{\tilde f_n})$ equals 
$g$ almost everywhere, and is hence $\mu$-measurable; of course all such are 
point-wise cluster points of $(\hat f'_n)$, a contradiction.
\end{proof}

\begin{corC}
$\ell^1$ does not embed in $X$ iff there exists a weakly precompact subset 
of $X$ isomorphically norming $X^*$.
\end{corC}

\begin{proof} 
If $\ell^1$ does not embed in $X$, $W = \Ba X$ is weakly precompact, by 
the $\ell^1$-Theorem.
Suppose conversely $W$ is a weakly precompact subset of $X$, isomorphically 
norming $X^*$.
It is easily seen that $\tilde W \defeq \{\alpha w:|\alpha|=1\}$ 
is also weakly pre-compact and hence $\overline{\co} \tilde W$,
the closed convex hull of $\tilde W$, is weakly 
pre-compact by Corollary~B. 
The Hahn-Banach Theorem (geometrical form) implies $\overline{\co}\tilde W$ 
contains $E \Ba X$ for some $\ep>0$, showing that every bounded sequence 
in $X$ has a weak-Cauchy subsequence.
\end{proof}

It is interesting to compare the equivalences 1--3  of Theorem~1
with the following 
equivalences, given in \cite{DU} 
(see Theorem~8, page~175 of \cite{DU}; $3\To 1$ is due to Uhl and 
$1\To 3$ is due to Stegall).

\begin{Thm}		
Let $X$ be given. 
Then the following are equivalent.
\begin{itemize}
\item[1.] $X^*$ has the Radon-Nikodym property (RNP).
\item[2.] Every integral operator from $Y$ to $X^*$ is nuclear, for all $Y$.
\item[3.] $Y^*$ is separable for all separable $Y\subset X^*$.
\end{itemize}
\end{Thm}

However this does not yield the sharp equivalence analogous to the one 
we obtain in $1\Leftrightarrow 3$ in Theorem~\ref{thm:main}.

\begin{prob}
Let $X$ be a separable Banach space. 
When does there exist an integral non-nuclear operator on $X^*$?
\end{prob}

\noindent Now if $X^*$ has the Radon-Nikodym  property, 
then by the above Theorem, every operator on $X^*$ is nuclear.
But also, if $X^{**}$ has the RNP, then: $T:X^*\to X^*$ integral implies 
$T^* :X^{**} \to X^{**}$ integral  implies $T^*$ 
is nuclear implies $T$ is nuclear, since $X^*$ is contractively complemented in 
$X^{***}$. 
Now if $X^*$ has the approximation property (ap) in addition, then (again by 
Grothendieck's results), 
$(K(X))^* = (X^* \mathop{\oplus}\limits^\vee X )^* = 
X^{**} \mathop{\oplus}\limits^\wedge X^* $ isometrically. 
This is so because Grothendieck's results yield that $(K(X))^* = I(X^*)$, 
the space of integral operators on $X^*$ (because $X$ has the ap). 
But then in the case $X^{**}$ has the RNP, we also have that the integral 
norm of $T$ on $X^*$  equals the integral norm of $T^{**}$, which equals
the nuclear norm of $T^{**}$, which equals the nuclear norm of $T^*$, 
which is then the same as its projective tensor norm. 
Also, $(X^{**} \mathop{\oplus}\limits^\vee X^*)^* = \L (X^{**})$ 
isometrically. 
We may thus summarize this discussion as follows.

\begin{prop}\label{prop3}
Let $X$ be a Banach space such that $X^*$ has the ap and $X^*$ or $X^{**}$ 
have the RNP. 
Then every integral  operator on $X^*$  is nuclear, $K(X)^* = X^{**}
\mathop{\oplus}\limits^\wedge X^*$ isometrically and 
$K(X)^{**} = \L (X^{**})$ isometrically.
\end{prop}

We show later on, however, that there exist separable Banach space $X$ 
so that $X^{**}$ has the metric approximation property, $X$ and $X^*$ 
fail the RNP, yet still every integral operator on $X^*$ is nuclear.

\begin{Rem}
The above discussion also shows that {\em given $X$ and $Y$, then if 
$X^*$ or $Y^*$ has the ap and either $X^*$ or $Y^*$ has the RNP, then 
$(X \mathop{\oplus}\limits^\vee Y)^* = X^* \mathop{\oplus}\limits^\wedge Y^*$
isometrically}. 
A slightly weaker result then this is due to Grothendieck
(\cite{Gr2}; see also \cite{DFS}).
\end{Rem}

We now pass to a detailed discussion of the following problem:

\begin{prob}
Under what conditions on Banach spaces $X$ and $Y$ is it so that $\ell^1$ 
embeds in $X\vtimes Y$, the injective tensor product of $X$ and $Y$?
\end{prob}

We solve this problem, for separable Banach spaces $X$ and $Y$ such that 
$X^*$ or $Y^*$ have the bounded approximation property (bap), in 
Theorem~\ref{thm:bap} below. 
This also yields a partial answer to Problem~2. 
(Of course under the ap assumption, the problem reduces to the study of 
separable spaces anyway, because if e.g. $X^*$ has the bap,
then we show in Lemma~9 that given $X_0$ a separable subspace of $X$,
there exists a separable subspace $X_1$, of $X$ with $X_1\supset X_0$ such 
that $X_1^*$ has the bounded approximation property. 
Thus if $\ell^1$ 
embeds in $X\vtimes Y$ and $X^*$ or $Y^*$ have the bap, there must exist 
$X_1,Y_1$ separable subspaces of $X$  and $Y$ so that $X_1^*$ 
or $Y_1^*$ has the bap and $\ell^1$ embeds in $X_1\vtimes Y_1$).

We shall freely use the standard (but rather non-trivial!) results
concerning tensor products of Banach spaces, due to Grothendieck \cite{Gr2},
as also exposed in \cite{DU}. 
For Banach spaces $X$ and $Y$, $I(X,Y)$ denotes the space of integral 
operators from $X$ to $Y$ and $N(X,Y)$ the space of nuclear operators from 
$X$ to $Y$.
For the definitions, including the norms, of these spaces, see the above 
references. 

Our next result provides as with sufficient conditions for $\ell^1$ 
to {\em not} embed in $X\vtimes Y$.

\begin{cor}\label{cor:Lewis}
Let $X$ and $Y$ be given Banach spaces, neither containing an 
isomorph of $\ell^1$, such that $X^*$ or $Y^*$ has the RNP. 
Then $\ell^1$ does not embed into $X\vtimes Y$.
\end{cor}

\begin{proof}
It obviously suffices to prove this for the case where $X$ and $Y$ are 
both separable.
By a result of Grothendieck, we also have that 
\begin{equation}\label{eq63}
(X\vtimes Y)^* = I(X,Y^*) = N(X,Y^*)
\end{equation}
(via trace duality). 
Furthermore, 
\begin{equation}\label{eq64}
N(X,Y^*)\ \text{ is isometric to a quotient space of }\ X^* \Atimes Y^*
\end{equation}
which implies that 
\begin{equation}\label{eq65}
[N(X,Y^*)]^*\ \text{ is isometric to a subspace of }\ 
X^* \Atimes Y^* = \L (X^*, Y^{**})\ .
\end{equation}
Now assuming $Y^*$ has the RNP, then $Y^*$ is separable, (by \cite{S}).
But since $Y^*$ is a separable Banach space, $Y^{**}$ is isometric to a 
subspace of $\ell^\infty$, and hence by Theorem~\ref{thm:main}, part~11, 
$\L(X^*,Y^{**})$ has cardinality $\G$, and so we have proved that 
$(X\vtimes Y)^{**}$ has cardinality $\G$, using 
\eqref{eq63}--\eqref{eq65}. 
Thus $\ell^1$ does not embed in $X\vtimes Y$ by Theorem~\ref{thm:main}, 
part~15 (i.e., by \cite{OR}).
\end{proof}

We shall show later on that Corollary~\ref{cor:Lewis} does not solve 
problem~2. 
In fact, the following consequence  of the result of R.~Haydon mentioned 
above is crucial for the solution.

\begin{prop}\label{prop:Haydon} 
Let $X$ and $Y$ be Banach spaces such that $\ell^1$ does not embed in 
$X\vtimes Y$. 
Then $I(X,Y^*)$ equals the closure of the finite rank operators from $X$ 
to $Y^*$ (endowed with the integral norm on $I(X,Y^*)$.
\end{prop}

\begin{proof}
Let $K = \{x^* \otimes y^* : \|x^*\|$, $\| y^*\| \le 1$, 
$x^* \in X^*$, $y^*\in Y^*\}$.
Then $K$ is a $w^*$-compact subset of $\Ba (X\vtimes Y)^*$ which 
isometrically norms $X\vtimes Y$. 
By the argument given in the proof of $1\To 9\To 10$ of 
Theorem~\ref{thm:main}, it follows from \cite{Hay} that 
$[K] = (X\vtimes Y)^* = I(X,Y^*)$. 
This proves Proposition~\ref{prop:Haydon}.
\end{proof}

\begin{remarks}
1. For $X$ and $Y$ separable, Proposition~\ref{prop:Haydon} follows 
by Theorem~\ref{thm:main}, part~10. 

2. Of course this result holds if we interchange $X$ and $Y$ in its 
statement, which is more naturally given in the language of tensor products:
{\em For separable Banach spaces $X$ and $Y$, if $\ell^1$ does not embed 
in $X\vtimes Y$, then $(X\vtimes Y)^*$ equals the closure of $X\otimes Y$ 
in the space of integral bilinear forms on $X^* \times Y^*$.}

3. It is apparently an open problem if the nuclear and integral norms 
coincide or are equivalent on $F(X,Y^*)$, $X$, $Y$ given 
($F(Z,W)$ denotes the space of finite rank operators from $Z$ to $W$).
Of course if this should be the case, for $X,Y$ satisfying the hypotheses 
of the Proposition, then its conclusion can be strengthened to state:
{\em Then every integral operator from $X$ to $Y^*$ is nuclear.}
\end{remarks}

We may now give a  definitive solution to Problem~2, under an 
approximation property assumption.

\begin{thm}\label{thm:bap} 
Let $X$ and $Y$ be separable Banach spaces such that $X^*$ or $Y^*$ has
the bounded approximation property. 
Then the following are equivalent
\begin{itemize}
\item[1.] $\ell^1$ does not embed in $X\vtimes Y$.
\item[2.] $\card \L(X^*,Y^{**}) = \G$.
\item[3.] $\card \L(X^*,Y^{**}) < 2^\G$.
\end{itemize}
Moreover when this occurs, every integral operator from $X$ to $Y^*$ 
is nuclear, and consequently 
$$(X\vtimes Y)^*  = X^* \Atimes Y^* \ ,\ \text{ hence }\ 
(X\vtimes Y)^{**} = \L(X^*, Y^{**})\ .$$
\end{thm} 

\begin{Rem}
As in the preceding result, we may (obviously) interchange $X$ and $Y$ 
in the statement of Theorem~\ref{thm:bap}.
\end{Rem}

\begin{proof}
Suppose first that 1 holds. 
The approximation property assumption insures that the integral and 
nuclear norms are equivalent on $F(X,Y^*)$, and these in turn are 
equivalent to the projective tensor product norm.
The final statement now follows from Proposition~\ref{prop:Haydon} and
the fact (due to Grothendieck) that $(X\vtimes Y)^* = I(X,Y^*)$; of course 
then 2 follows by Theorem~\ref{thm:main} (i.e., by \cite{OR}).

$2\Rightarrow 3$ is trivial.

Now suppose that 3 holds, and assume that
\begin{equation}\label{eq66} 
X^* \ \text{ has the approximation property .}
\end{equation}
Let $K = \{ x^* \otimes y^* : (x^*, y^*) \in \Ba X^* \times\Ba Y^*)\}$, 
endowed with its $w^*$-topology as a subset of $I(X,Y^*) = (X \vtimes Y)^*$.
Thus $K$ is a $w^*$-compact isometrically norming subset of $(X\vtimes Y)^*$.
To prove that 1 holds it suffices 
to prove that given $(A_n)$ a bounded sequence in $X\vtimes Y$, 
then defining $\hat A_n(x^* \otimes y^*) = \langle x^* \otimes y^*,A_n\rangle$
for all $n$ and $x^* \otimes y^* \in K$, then 
\begin{equation}\label{eq67}
\hat A_n \ \text{ has a pointwise convergent subsequence.}
\end{equation}
If \eqref{eq67}  is false, then by 19 of Theorem~\ref{thm:main}, 
\begin{equation}\label{eq68} 
\text{The family $\F$ of pointwise cluster points of $(\hat A_n)$ 
has cardinality $2^\G$.}
\end{equation}
(Of course we can directly apply Lemma~2 to also see that \eqref{eq68} 
holds.)
Now suppose that $f:K\to \KK$ is a pointwise cluster point of $(\hat A_n)$ 
(where $\KK$ denotes the scalar field), and suppose that $\|A_n\| \le C$
for all $n$, $C<\infty$. 
Now we identify $X\vtimes Y$ with $K_* (X^*,Y)$, the Banach space of 
compact operators in $\L (X^*,Y)$ which are weak*-norm continuous on 
bounded subsets of $X^*$. 
(This is legitimate by results in \cite{Gr2}.) 
Of course trivially $K_* (X^*,Y)\subset \L(X^*,Y^{**})$. 
Then we claim that 
\begin{eqnarray}
&&\text{there exists a unique $T_f\in \L(X^*,Y^{**})$ such that}\label{eq69}\\
&&\text{(i)}\ \qquad \langle T_f(x^*),y^*\rangle = f(x^*\otimes y^*)
\ \text{ for all }\ (x^*,y^*) \in \Ba X^*\times \Ba Y^* \notag\\
\noalign{\text{and moreover}}
&&\text{(ii)}\qquad  
T_f \ne T_g\ \text{ if }\ f,g\in \F\ ,\qquad f\ne g\ .  \notag
\end{eqnarray}
Indeed, we may choose a net $(n_\alpha)_{\alpha \in D}$ so that 
\begin{equation}\label{eq70} 
\lim_\alpha A_{n_\alpha} (x^* \otimes y^*) = f(x^*\otimes y^*) 
\ \text{ for all } \  x^* \otimes y^* \in K\ .
\end{equation}
But then it follows easily that 
\begin{equation}\label{eq71} 
\lim_\alpha A_{n_\alpha} (x^*\otimes y^*) \defeq G_f (x^*,y^*) 
\ \text{ exists for all }\ 
(x^*,y^*) \in X^* \times Y^*\ ,
\end{equation}
and moreover
\begin{equation}\label{eq72}
\text{$G_f$ is a bilinear form on $X^*\times Y^*$ with norm bounded by $C$.}
\end{equation} 
Of course then there is a unique $T_f \in \L(X^*,Y^{**})$ such that 
\begin{equation}\label{eq73} 
\langle T_f (x^*),y^*\rangle = G_f (x^*,y^*)\ \text{ for all }\ 
(x^*,y^*) \in X^* \times Y^*
\end{equation}
and moreover it is obvious that $f\ne g\in \F\Rightarrow G_f\ne G_g 
\Rightarrow T_f \ne T_g$.
Thus $\card \L(X^*,Y^{**}) = 2^\G$, contradicting 3.

The final statement of Theorem~\ref{thm:bap} follows by 
Proposition~\ref{prop:Haydon}, Remark~3 following its proof, and the 
theorem of Grothendieck that $X^* \Atimes X^* = N(X,Y^*)$ since $X^*$ 
has the bap.
\end{proof}

\begin{cor}\label{cor7}
Let $X$ and $Y$ satisfy the hypotheses of Theorem~\ref{thm:bap} and 
suppose $\ell^1$ does not embed in either space. 
Then if every operator from $X^*$ to $Y^{**}$ has separable range, $\ell^1$ 
does not embed in $X \vtimes Y$. 
In particular, this holds if every operator from $X^*$ to $Y^{**}$ is compact.
\end{cor}

\begin{proof} 
Since $\ell^1$ does not embed in $Y^{**}$, $\card Y^{**} = \G$, and hence 
the cardinality of the family of countable subsets of $Y^{**}$ is $\G$.
If $Z$ is a separable subspace of $Y^{**}$, we may 
choose a countable dense subset $Z_0$ of $Z$, and so we have established that 
\begin{equation}\label{eq74}
\text{the cardinality of the family of separable subspaces of $Y^{**}$ 
equals $\G$.}
\end{equation}
Now for each separable non-zero subspace $Z$ of $Y^{**}$, 
\begin{equation}\label{eq75}
\L(X^*,Z)\ \text{ has cardinality $\G$}
\end{equation}
by condition 13 of Theorem~\ref{thm:main}.
Since the union of a family of sets of cardinality $\G$ has cardinality $\G$ 
provided each of the sets has cardinality $\G$, 
\eqref{eq74} and \eqref{eq75} imply that 3 of Theorem~\ref{thm:bap} holds, 
and hence $\ell^1$ does not embed in $X \vtimes Y$ by Theorem~\ref{thm:bap}.
The final statement of the Corollary now trivially follows.
\end{proof}

\begin{prob}\label{prob3}
Let $X$ and $Y$ satisfy the hypotheses of Theorem~\ref{thm:bap}, and 
suppose every integral operator from $X$ to $Y^*$ is nuclear. 
Is it so that $\ell^1$ does not embed in $X\vtimes Y$ unless $\ell^1$ embeds 
in $X$ or $Y$?
\end{prob}

Problem~\ref{prob3}, in turn, suggests 

\begin{prob}\label{prob4}
If $X$ and $Y$ are given such that $L^1$ embeds in $X^* \Atimes Y^*$; does 
$L^1$ embed in $X^*$ or $Y^*$?
\end{prob}

By the results in \cite{H1} and \cite{P}, an affirmative answer holds 
for particular $X,Y$ iff 
$\ell^1$ embeds in $X$ or $Y$.
An affirmative answer to Problem~\ref{prob4} implies 
an affirmative answer to Problem~\ref{prob3}. 
For suppose to the contrary, that $X,Y$ provide a counterexample to 
Problem~\ref{prob3}. 
Since $\ell^1$ embeds in $X\vtimes Y$, $L^1$ embeds in 
$(X\vtimes Y)^* = I(X,Y^*)$. 
But since $\ell^1$ does not embed in $X$ or $Y$, $L^1$ does not embed in 
$X^*$ or $Y^*$, and so if the answer to the second part of \ref{prob4} 
is affirmative, $L^1$ does not embed in $X^* \Atimes Y^*$, which thus cannot 
equal $I(X,Y^*)$ contradicting the final statement in Theorem~\ref{thm:bap}. 
(See remark~4 following proof of Theorem~\ref{thm11} for further comments 
on Problem~\ref{prob4}.)  

\begin{Rem}
If we do not deal with dual spaces in Problem~\ref{prob4}, then the 
answer is negative.
A remarkable result of Talagrand \cite{T} asserts the existence of separable
spaces $X$ and $Y$ so that $L^1$ embeds in $X\oplus Y$ but $L^1$ does 
not embed in $X$ or $Y$.
It follows, since $L^1$ is isomorphic to a finite codimensional subspace
of itself, 
that also $L^1$ embeds in $X\oplus Y_0$ if $Y_0$ is of codimension~1 in $Y$.
Let $x_0,y_0$ be norm one elements of $X^*,Y^*$ respectively; choose 
$x_0^*,y_0^*$ norm one elements of $X^*,Y^*$ with $x_0^* (x_0) = y_0^*(y_0)=1$,
and let $\tilde Y = y_0^\perp$. 
Then $X\otimes [y_0]$ and $[x_0] \otimes \tilde Y$ are isometric to $X$ 
and $\tilde Y$ in $X\Atimes Y$, and $I\otimes(y_0^* \otimes y_0)$ is a 
contractive projection from $X\Atimes Y$ onto  $X\otimes [y_0]$ with 
$[x_0] \otimes \tilde Y$ contained in its kernel. 
It follows that $(X\oplus [y_0]) + ([x_0]\otimes \tilde Y)$ is isomorphic to 
$X\oplus \tilde Y$ in $X\Atimes Y$ and  thus $L^1$ embeds in 
$X\Atimes Y$.
\end{Rem}

Of course this counterexample cannot simply lead to a counterexample for 
Problem~4, because if $L^1$ embeds in $X^*\oplus Y^*$, $\ell^1$ embeds in 
$X\oplus Y$ by \cite{H1} and \cite{P}, and then obviously $\ell^1$ embeds in 
$X$ or $Y$.
Nevertheless, I am inclined to believe that the answer to problem~3 and hence
to problem~4 is negative. 
Leter on, we give examples where in fact $\ell^1$ does not embed in $X^*$ 
or $Y^*$, both $X^*$ and $Y^*$ have the metric approximation property and
fail the RNP, and indeed $(X\vtimes Y)^* \ne X^*\Atimes Y^*$.
However the proof is rather delicate, using specific properties of these 
spaces, and I'm inclined to believe there is no technique general enough 
to give an affirmative answer to problem~4.

The next result summarizes consequences of the previous results in the 
context of Problem~2.  

\begin{thm}\label{thm8}
Let $X$ be a given Banach space.
Consider the following two properties
\begin{itemize}
\item[{}]  $P_1$. Every integral operator on  $X^*$ is nuclear.
\item[{}]  $P_2$. $\ell^1$ does not embed in $K(X)$.
\item[1.] If $P_1$ holds, $\ell^1$ does not embed in $X$.
\item[2.] If $X^*$ has the approximation property and $X^*$ or $X^{**}$ has
the RNP, then $P_1$ holds.
\item[3.] If $X^*$ has the bounded approximation property and $P_2$ holds, 
$P_1$ holds.
\item[4.] If $X^*$ has the bounded approximation property, then $P_2$ holds 
if and only if for all separable subspaces $Z$ and $Y$ of $X^*$ and $X$ 
respectively such that $Z^*$ has the bounded  approximation property, 
the cardinality of $\L (Z^*,Y^{**})$ is less than $2^\G$ (iff the 
cardinality equals  $\G$).
\end{itemize}
\end{thm}

1. This follows from Theorem~\ref{thm:main}, part~4, for if 1 is false,
this implies there exists an integral operator on $X^*$ which is not 
compact, hence not nuclear.

2. This is Proposition~\ref{prop3}.

3. If $X^*$ has the bap, then $K(X) = X^*\vtimes X$ (because $X$ has the ap), 
and hence 
$$K(X)^* = I(X^*) = X^{**} \Atimes X^*$$
by Proposition~\ref{prop:Haydon}, since the hypotheses imply that the 
integral, nuclear, and projective tensor norms are equivalent on 
$X^{**} \otimes X^*$. 
Of course the final equality implies that $P_1$ holds.

To show 4, we need 

\begin{lem}\label{lem9}
If $X^*$ has the bap, then very separable subspace $E$ of $X$ is contained 
in a separable subspace $Z$ of $X$ such that $Z^*$ is isomorphic to a 
complemented subspace of $X^*$.
\end{lem} 

We first complete the proof of Theorem~\ref{thm8}, then prove Lemma~\ref{lem9}.

1.  As in the proof of 3, $K(X) = X^*\vtimes X$ and so the ``only if'' 
statement follows immediately from Theorem~\ref{thm:bap}. 
To see ``if'', suppose to the contrary that $\ell^1$ embeds in $K(X)$. 
Then we can choose separable subspaces $E$ and $Y$ of $X$ and $X^*$ 
respectively such that $\ell^1$ embeds in $E\vtimes Y$. 
But if we choose $Z$ satisfying the conclusion of Lemma~\ref{lem9}, then 
$Z^*$ has the bap, and then $\ell^1$ does not embed in $Z\vtimes Y$, by 
Theorem~\ref{thm:bap}, a contradiction.\qed
\medskip

\begin{proof}[Proof of Lemma~\ref{lem9}]
Since $X^*$ has the bap, so does $X$, and so we may choose $1<\lambda<\infty$
so that for all finite-dimensional subspaces $F$ of $X$, there exists
a finite rank operator $T$ on $X$ such that
\begin{equation} 
\|T\| \le \lambda\ \text{ and }\ T\mid F =  I\mid F\ .
\label{tag-a}
\end{equation}
Let $e_0,e_1,e_2,\ldots$ be an enumeration of a countable dense subset of $E$.
Choose $T_1$ a finite rank operator on $X$ satisfying \eqref{tag-a}, 
where $F= [e_0,e_1]$ and $\text{``$T$''} = T_1$. 
Suppose $n\ge 1$ and finite rank operators $T_1,T_n$ have been chosen. 
Now let $F = [\text{Range }T_n,[e_{n+1}]]$. 
Choose $T_{n+1}$ satisfying \eqref{tag-a} with $\text{``$T$''} = T_{n+1}$.

This completes the inductive construction of the $T_n$'s; 
Let $F_n = T_n (X)$ and $Z = \overline{\cup F_j}$. 
Note that $F_n\subset F_{n+1}$ and $e_n\in F_n$ for all $n$; 
hence $E\subset Z$. 
It follows that 
\begin{equation}
\text{$T_n(X)\subset Z$ for all $n$ and $T_n (z) \to z$ in norm for all 
$z\in Z$\ .}	
\label{tag-b}
\end{equation}
Now the compactness of the unit ball of $X^*$ in the weak* topology implies 
that there exists a net $(T_{n_\alpha})_{\alpha\in\infty}$ 
(i.e., a subnet  of the sequence $(T_j)$ and an operator $P$ on $X^*$ 
such that
\begin{equation}
T_{n_\alpha}^* (x^*) \to P(x^*)\  w^* \ \text{ for all }\ x^* \in X^*\ .
\label{tag-c}
\end{equation}
\eqref{tag-c} shows that $\|P\| \le \lambda$ (so $P$ is indeed bounded).
Now if $x^* \in Z^*$, then \eqref{tag-b} implies that $T_n^* (x^*) =0$ 
for all $n$, which implies 
\begin{equation}
P(x^*) =0\ \text{ for all }\ x^* \in Z^+\ .	
\label{tag-d}
\end{equation}

Now by  our construction, 
$T_{n+1} \mid\text{Range }T_n = I\mid\text{Range }T_n$ for all $n$, which 
implies that for all $k$ and $n> k$, $T_n T_k = T_k$, and so taking 
adjoints, $T_k^* T_n^* = T_k^*$ for all $k$ and $n> k$. 
But then for $x^* \in X^*$, we deduce, thanks 
to the weak*-continuity of $T_k^*$, 
that 
\begin{equation}
\lim_\alpha T_k^* T_{n_\alpha}^* (x^*) = T_k P(x^*)  = T_k (x^*)\ .  
\label{tag-e}
\end{equation}
In turn, \eqref{tag-e} implies, after taking another limit, that $P^2 =P$; 
hence $P$ is a projection. 
Finally, 
\begin{equation}
\text{if $P(x^*) =0$, then $x^* \in Z^\perp$.}	
\label{tag-f}
\end{equation}
For if not, we would find $z\in Z$, $\langle x^*,z\rangle \ne0$. 
But then 
\begin{equation}
\langle Px^*,z\rangle = \lim_\alpha \langle T_{n_\alpha}^*,x^*,z\rangle 
= \lim_\alpha \langle x^*,T_{n_\alpha}z\rangle = (x^*,z)\ne0\ .  
\label{tag-g}
\end{equation}
\eqref{tag-d} and \eqref{tag-f} show that $Z^\perp$ is the kernel of $P^*$.
Thus it follows that $PX^*$ is isomorphic to $X^*/Z^\perp$, which of course 
is isometric to $Z^*$.
\end{proof}

\begin{Rem}
It is known that the conclusion of the Lemma holds without any approximation
property assumptions, but we included a proof of our needed result for 
completeness.
We only used the hypothesis of the Lemma to obtain that $X$ has the bounded 
approximation property.
With a little more care in the proof, we may thus obtain the following result:
{\em If $X$ has the $\lambda$-bap and $E$ is a separable subspace of $X$, 
there is a separable subspace $Z$ of $X$ with the $\lambda$-bap, containing 
$E$, such that $Z^*$ is $\lambda$-isomorphic to a $\lambda$-complemented 
subspace of $X^*$.}
\end{Rem}

Of course Theorem~\ref{thm8} suggests a possible solution to Problem~2, 
namely that $P_1$ and $P_2$ are equivalent properties, if $X^*$ has the bap.

\begin{prob}\label{prob5} 
If $X^*$ has the bap, does $P_1$ imply $P_2$?
\end{prob}

Of course Problems 4 and 6 are linked, for an affirmative answer to Problem~4
implies an affirmative answer to Problem~6. 

We now give examples illustrating Theorem~\ref{thm:bap}. 
The examples show that the conditions of Corollary~\ref{cor:Lewis} 
are not necessary, to insure that $\ell^1$ does not embed in $X\vtimes Y$. 
The examples also show that there exist separable Banach spaces $X$ and $Y$
not containing $\ell^1$ isomorphically, so that there exists an integral 
non-nuclear operator from $X$ to $Y^*$. 
The examples are natural generalizations of the James tree space $JT$.
We first define these spaces and summarize some of their properties.

Let $\D$ be the dyadic tree; that is, $\D$ is the set of all finite 
sequences of 0's and 1's, ordered by extension. 
Let $\tau :\D\to N$ be the standard bijection, given as 
\begin{equation}\label{eq76}
\emptyset,\ (0),\ (1),\ (00),\ (01),\ (10),\ (11),\ \cdots
\end{equation}
and let $(e_j)$ be the unit vectors bases of $c_{0,0}$ the space of all 
sequences of scalars which are ultimately zero. 
A non-empty subset $s$ of $\D$ is called a {\em segment} provided 
\begin{equation}\label{eq77} 
\text{$s$ is totally ordered, and whenever $\alpha<\gamma<\beta$ in 
$\D$ with $\alpha,\beta\in s$, then $\gamma\in s$.}
\end{equation}
A maximal totally ordered subset of $\D$ is called a {\em branch} 
of $\D$; equivalently, this is a segment containing $\emptyset$, 
unbounded above; branches can also be identified with infinite sequences of 
0's and 1's, where if $(\ep_n)$ is such a sequence, the corresponding 
branch is $\{\emptyset, (\ep_1), (\ep_1,\ep_2), (\ep_1,\ep_2,\ep_3),\ldots\}$. 

Given a segment $s$ of $\D$, define a functional $s^*$ on $c_{00}$ by 
\begin{equation}\label{eq78}
s^* (f) = \sum_{j\in s} f_{\tau (j)}\ .
\end{equation}
Now fix $p$, $1<p<\infty$, and define a norm $\|\cdot\|_{JT_p}$ on $c_{00}$ by 
\begin{equation}\label{eq79} 
\|f\|_{JT_p} = \max \bigg\{ \bigg( \sum_{i=1}^k |s_i^* (f) |^p\bigg)^? 
: k\ge 1
\text{ and }s_1,\ldots,s_k\text{ are disjoint segments of }\D\bigg\}\ .
\end{equation}
We define $JT_p$ to be the completion of $(c_{00},JT_p)$. 
$JT_2$ is the space defined as $JT$ in \cite{LS}, 
discovered in \cite{J3}, and the proofs of the 
properties of $JT_p$ that we require are the same as those for $JT$.
(For an alternate exposition, see \cite{B}.)

Let us first note that 
\begin{equation}\label{eq80}
\text{$s^*$ is a norm one functional on $JT_p$ for any segment $s$.}
\end{equation} 
We refer to the functionals $\beta^*$, for $\beta$ a branch of $\D$, as 
{\em branch functionals}.
For $W$ a non-empty subset of $\D$, set 
\begin{equation}\label{eq81} 
\widetilde W = [e_j : \tau (j) \in W]\ .
\end{equation}
Then it is obvious that for any segment $s$, 
\begin{equation}\label{eq82} 
\tilde s\text{ is contractively complemented in $JT_p$, with complement 
$\widetilde{\D\sim s}$.}
\end{equation}
Let us also note that given a branch $\beta$ of $\D$, 
$\beta\sim (\ep_j)_{j=1}^\infty$, then if we let 
$\beta_\alpha =  e_{\tau^{-1} (\ep_1,\ldots,\ep_{n-1})}$, $j=1,\ldots$,
then 
\begin{equation}\label{eq83}
\begin{split}
&(\beta_j)_{j=1}^\infty \ \text{ is isometrically equivalent to the 
boundedly complete basis for } J\\
&\text{where $J$ denotes the quasi-reflexive order one space of James.}
\end{split}
\end{equation}
In particular, $(\beta_j)_{j=1}^\infty$ dominates the summing basis, and the 
branch functional $\beta^*$ restricted to $\tilde \beta$ is the summing 
functional. 

We denote by $\D_n$ the family of elements of $\D$ of length $n$; 
then letting $\alpha_1^n,\ldots,\alpha_{2^n}^n$ be an enumeration of $\D_n$, 
and setting $e_j^n = e_{\tau^{-1} (\alpha_j^n)}$ for all $j$, we have that 
\begin{eqnarray}
&&\tilde\D_n\ \text{ is contractively complemented in $JT_p$, by 
$\widetilde{\sim \D_n}$\ ,}\label{eq84}\\
&&\text{and $(e_j^n)_{j=1}^{2^n}$ is isometrically equivalent to the 
$\ell_{2^n}^p$ basis.}\notag
\end{eqnarray}
We now summarize the rest of the properties of $JT_p$ that we shall need 
(some of which are decidedly non-trivial). 
(Throughout, $r^*$ denotes the conjugate index to $r$, 
$\frac1r +\frac1{r^*}=1$.) 

\begin{thm}\label{thm10}
$\quad$
\begin{itemize}
\item[1.] $(e_j)$ is a boundedly complete monotone basis for $JT_p$.
\item[2.] $(JT_p)^*$ is the norm-closed linear space of the functionals 
$(e_j^*)$ biorthogonal to $(e_j)$, and the branch functionals.
\item[3.] The set of branch functionals $\B$ is homeomorphic  to the 
Cantor set, in the weak*-topology. 
If for $\alpha \in\D$, $\B_\alpha = \{\beta^* :\alpha \in\beta\}$, 
then $(\B_a)_{\alpha \in\infty}$ is a family of $w^*$-clopen subsets 
of $\B$, forming a base for its topology.
\item[4.] Letting $(JT_p)_*$ be the closed linear span of the biorthogonal 
functionals $(e_\alpha^*)$, then 
\begin{equation}\label{eq85}
(JT_p)^* /(JT_p)_*\ \text{ is isometric to }\  \ell_\G^{p^*}\ .
\end{equation}
\item[5.] Every normalized weakly null sequence in $JT_p$ has a subsequence 
equivalent to the $\ell^p$-basis.
\end{itemize}
\end{thm}

(5. was established for $JT$ in \cite{AI}; see \cite{B} for an alternate 
treatment of the proof of Theorem~10 for $JT$, which immediately generalizes 
to the case of $JT_p$.)

We may draw some immediate  consequences. 
First, since $(JT_p)_*$ is thus canonically embedded in its double dual, 
$JT_p$ is canonically complemented in $(JT_p)^*$, by $((JT_p)_*)^\perp$. 
Hence it follows from Theorem~10		
part~4 that 
\begin{equation}\label{eq86} 
JT_p^{**} \ \text{ is isometric to }\ JT_p \oplus \ell_\G^p\ .
\end{equation}
Of course $JT_p$ has the Radon-Nikodym property, being a separable dual.
Then obviously by \eqref{eq86}, $JT_p^{**}$ also has the Radon-Nikodym 
property and the metric approximation property. 
So we deduce that 
\begin{equation}\label{eq87}
(JT_p)_*,\ JT_p,\ (JT_p)^*,\ \text{ and }\ JT_p^{**}\ \text{ all have the 
metric approximation property.}
\end{equation}
Of course it also follows immediately that $\ell^1$ does not 
embed in $JT_p^{**}$.

We may now formulate our final main result; which shows in particular that 
for $1<p<\infty$, 
\begin{equation}\label{eq88}
\ell^1 \ \text{ embeds in }\ JT_p \vtimes JT_p\ \text{ iff }\ 2\le p<\infty\ .
\end{equation} 
(I am indebted to J.~Diestel for showing me the important special case: 
$\ell^1$ embeds in $JT\vtimes JT$.)

\begin{thm}\label{thm11} 
Let $1<p,q<\infty$. 
Then the following are equivalent.
\begin{itemize}
\item[1.] $\ell^1$ embeds in $JT_p \vtimes JT_q$.
\item[2.] $p^* \le q$.
\item[3.] There exists an integral non-nuclear operator from $JT_p$ to 
$(JT_q)^*$.
\end{itemize}
\end{thm} 

\begin{Rem}
Of course since $JT_p$, $JT_q$ satisfy the hypotheses of Theorem~\ref{thm:bap},
we may add as a 4th equivalence
\begin{itemize}
\item[4.]  {\em $\card \L_{p,q} = 2^\G$, where $\L_{p,q} = \L((JT_p)^*,
(JT_p)^{**})$.}
\end{itemize}
\end{Rem}

\begin{proof} 
Suppose first that $p^* \le q$. 
We show that 3 holds, which implies that 1 holds by Theorem~\ref{thm:bap}.
Let $\mu$ be the ``Cantor probability measure'' on the Borel sets of 
$\B^q$, associated to the family of clopen sets $(\B_\alpha^q)_{\alpha\in\D}$.
That is, $\mu$ is the unique Borel measure so that 
\begin{equation}\label{eq89}
\mu (\B_\alpha^q) = \frac1{2^n} \ \text{ if }\ \alpha \in \D_n
\ \text{ for }\ n=0,1,2,\ldots\ .
\end{equation} 
(Since we are dealing with different spaces, we are denoting by $\B^q$ the 
set of Banach functionals on $JT_q$, $\B_\alpha^q$ its associated 
clopen-set basis.)
Now we define a map $V:L^1(\mu) \to (JT_q)^*$ by 
\begin{equation}\label{eq90}
V f = \int_\B f(w) w^* \,d\mu(w)\ \text{ for }\ f\in L^1 (\mu)\ .
\end{equation}
The integral denotes the weak*-integral; thus, we easily verify that 
given $x\in JT_q$, then since $w^* \to w^* (x)$ is a continuous function 
of norm at most $\|x\|$ on $\B$, then for $f\in L^1 (\mu)$
\begin{equation}\label{eq91}
x\to \int_\B f(w) w^* (x)\,d\mu (w)
\end{equation} 
is in $(JT_q)^*$ and has norm at most $\|f\|_{L^1(\mu)}$ thus \eqref{eq90} 
is well defined, $V$ is indeed a linear operator, and in fact $\|V\|\le1$. 
Then we have that for all $\alpha \in \D$, 
\begin{equation}\label{eq92} 
\begin{cases}
(V \chi_{\B_\alpha}) (e_{\tau(\beta)}) = \frac1{2^{|\alpha|}}
&\text{if $\beta\in\B_\alpha$}\\
\noalign{\vskip6pt}
V\chi_{\B_\alpha} (e_{\tau (\beta)}) =0&\text{if $\beta\notin \B_\alpha$}
\end{cases}
\end{equation} 
where we set $|\alpha| =n$ if $\alpha \in \D_n$. 
Next,  let $\varphi :\B^q \to \B^p$ be the canonical  homeomorphism such that 
$\varphi (\B_\alpha^q) = \B_\alpha^p$ for all $\alpha$.
Now define $U:JT_p \to C(\B^q)$ by 
\begin{equation}\label{eq93}
(Ux) (\beta^*) = \varphi (\beta^*) (x)\ \text{ for all }\ x\in JT_p\ ,
\ \ \beta\in \B^q\ .
\end{equation}
\end{proof}

Then it is obvious that $U$ is a linear contraction and thus 
\begin{equation}\label{eq94} 
T = VU \ \text{ is an integral operator from $JT_p$ to }\ (JT_q)^*\ .
\end{equation}
Thus by trace duality, we have (by Grothendieck's fundamental theory) given 
\newline
$x= \sum_{i=1}^n x_i \otimes y_i \in JT_p\otimes JT_q$, then 
\begin{equation}\label{eq95} 
\Big| \sum_{i=1}^n Tx_i (y_i)\Big| \le 
\Big\| \sum_{i=1}^n x_i\otimes y_i\Big\|
\end{equation}
(and $\sum_{i=1}^n x_i \otimes  y_i\to \sum_{i=1}^n Tx_i (y_i) = 
\tr T\sum_{i=1}^n y_i\otimes x_i$ is a well defined 
linear functional; we used in \eqref{eq95}  that the integral norm of $T$
is at most~1). 
To prove that $T$ is integral but not nuclear, we shall show that  there 
exists a $G\in (JT_p\vtimes JT_q)^{**}$ such that 
\begin{equation}\label{eq96}
G(T) =1\ \text{ and }\ G(f) =0\ \text{ for all }\ f\in JT_p^*\Atimes JT_q^*\ .
\end{equation}
Of course to prove the last claim in \eqref{eq96}, it suffices to prove that 
\begin{equation}\label{eq97} 
G(x^* \otimes y^*) =0\ \text{ for all }\ x^* \in JT_p^*, \ y^*\in JT_q^*\ .
\end{equation}

Let us call $x^* \otimes y^*$ a {\em basic tensor}  if $x^*,y^*$ are each 
either a biorthogonal functional or a summing functional in $JT_p^*$, 
respectively $JT_q^*$. 
Thanks to part~2 of Theorem~10	
and the continuity of $G$, it actually 
suffices to prove that \eqref{eq97} holds for all basic tensors 
$x^*\otimes y^*$.

Now fix $n$, and define $A_n\in JT_p\otimes JT_q$ by 
\begin{equation}\label{eq98} 
A_n = \sum_{i=1}^{2^n} e_{\alpha_i^n} \otimes e_{\alpha_i^n}\ .
\end{equation} 
Now let $P_n$ be the canonical projection from $JT_p$ onto 
$[e_i^n]_{i=1}^{2^n}$; thus $P_n^*$ is the canonical projection from 
$(JT_p)^*$ onto $[e_i^{n*} ]_{i=1}^{2^n}$. 
Let $Q_n : [e_i^{n*}]_{i=1}^{2^n} \to [e_i^n]_{i=1}^{2^n}$ be the 
linear map such that $Q_n e_i^{n*}  = e_i^n$ for $1\le i\le 2^n$, where 
$e_i^{n*}$ denotes an element of $(JT_p)^*$ and $e_i^n$ denotes an element 
of $JT_q$. 
Then 
\begin{equation}\label{eq99} 
A_n = Q_n P_n\ .
\end{equation}
Thanks to \eqref{eq98}, $(e_i^{n*})_{i=1}^{2^n}$ is isometrically 
equivalent to the $\ell_{2^n}^{p^*}$ bases, $(e_i^n)_{i=1}^{2^n}$ is 
isometrically equivalent to the $\ell_{2^n}^q$ basis, and hence since 
2 holds, $\|Q_n\| =1$, $\|P_n\|=1$, and so by \eqref{eq99} 
\begin{equation}\label{eq100} 
\|A_n\| \le 1\ .
\end{equation}

We shall define $G$ satisfying \eqref{eq97} and \eqref{eq98} in two steps.
First, let $G_1 \in (JT_p\vtimes JT_q)^{**}$ be a $w^*$-cluster point of 
$(A_n)_{n=1}^\infty$. 
Of course \eqref{eq100} shows that $\|G_1\| \le 1$. 
We have that letting $\langle ,\rangle$ be the pairing between 
$I(JT_p,(JT_q)^*)$ and $JT_p \otimes JT_q$ given after (70) then 
for all $n$, 
\begin{eqnarray}
\langle T,A_n\rangle & = &\sum_{i=1}^{2^n} T(e_i^n) (e_i^n)\label{eq101}\\
&=& \sum_{i=1}^{2^n} \frac1{2^n}\ \text{ using (98)}\notag\\
&=& 1\ .\notag
\end{eqnarray}
Thus $G_1 (T)=1$.

Next, we claim that for all basic tensors $x^*$ and $y^*$, 
\begin{equation}\label{eq102}
G_1 (x^* \otimes y^*) =  \lim_{n\to\infty} \sum_{j=1}^{2^n} 
x^* (e_i^n) y^* (e_i^n)
\end{equation}
(part of this assertion is the claim that this limits exists) and so that 
setting $H=G_1$, then 
\begin{eqnarray}
&&H(x^*\otimes y^*) =0\text{ unless there is a branch $\gamma$ so that 
$x^* =\gamma^*$ (in $(JT_p)^*$),} \label{eq103}\\
&&\qquad y^* = \gamma^* \text{(in $(TJ_p)^*$), and then }
H (\gamma^* \otimes \gamma^*) =1\ .\notag
\end{eqnarray}
Now it is obvious that the limit in \eqref{eq102} is zero if one of $x^*$ or 
$y^*$ is a biorthogonal functional.
Suppose that there are branches $\gamma\ne\beta$ such that $x^*=\gamma^*$ 
and $y^* = \beta^*$, and choose $k$ so that the branches split at level 
$k$; i.e., if $\gamma = (\gamma_j)_{j=0}^\infty$ and  
$\beta = (\beta_j)_{j=0}^\infty$, then $\gamma_k \ne\beta_k$, which implies 
that $\gamma_n \ne \beta_n$ for all $n>k$. 
But then if $n>k$, there are unique $i\ne j$ such that $\alpha_i^n\in\gamma$,
$\alpha_j^n\in\beta$, and hence 
\begin{equation}\label{eq104}
\gamma^* (e_\ell^n) \beta^* (e_\ell^n) =0\ \text{ for all }\ 1\le \ell\le 2^n\ .
\end{equation}
Finally, if $x^* = \gamma^* = y^*$, then for each $n$, there is exactly one 
$i$ with $\alpha_i^n\in\gamma$, and then $\gamma^* (e_i^n)=1$, 
$\gamma^* (e_j^n) =0$, $j\ne i$, and hence $\sum_{\ell=1}^{2^n} \gamma^* 
(e_\ell^n) \gamma^* (e_\ell^n) =1$ for all $n$.
(Actually, $(A_n)$ converges in the $w^*$ operator topology on $\L_{p,q}$ 
to the operator $S$ such that $S(e_j^*)=0$ for all $j$, and $S(\gamma^*)= 
\gamma^{**}$ for all branches $\gamma$, where $\gamma^{**} (\beta^*)=
\delta_{\gamma\beta}$ all $\beta$ and $\gamma^{**} (e_j^*) = 0$ all $j$; 
$(\gamma^{**})_{\gamma\in\B}$ is in fact the basis for $\ell_\G^q$ in 
$(JT_q)^{**})$. 

Now let $\F$ be the family of all finite non-empty subsets of $\B$, directed
by inclusion. 
For each $n$ and $F= \{\gamma_1,\ldots,\gamma_n\}\in \F$ of cardinality $n$, 
choose 
\begin{equation}\label{eq105} 
m_n \ge n\ \text{ such that the branches $\gamma_1,\ldots,\gamma_n$ have 
split at level $m_n$}
\end{equation}
That is, we may choose $n$ different  integers $i_1,\ldots,i_n$ in 
$\{ 1,2,\ldots,2^{m_n}\}$ such that 
\begin{equation}\label{eq106}
\alpha_{i_j}^{m_n} \in \gamma_j \ \text{ for }\ 1\le j\le n\ .
\end{equation}
Define $B_F \in J_p \otimes J_q$ by 
\begin{equation}\label{eq107} 
B_F = \sum_{j=1}^n e_{i_j}^{m_n} \otimes e_{i_j}^{m_n}\ .
\end{equation}
We have, just as in the definition of the $A_n$'s, that 
\begin{equation}\label{eq108} 
\|B_F\| \le 1\ \text{ for all }\ F\in \F\ .
\end{equation}
Now let $G_2$ be a $w^*$-cluster point of $(B_F)_{F\in \F}$ in 
$(JT_p\otimes JT_q)^{**}$.
Then
\begin{equation}\label{eq109}
G_2 (T) =0\ .
\end{equation}
Indeed, we have that 
\begin{equation}\label{eq110}
\langle T,\beta_n\rangle  = \sum_{j=1}^n \frac1{2^{m_n}} \le 
\frac{n}{2^n} \to 0\ \text{ as }\ n\to\infty\ .
\end{equation}
Moreover, we claim that setting $H=G_2$, then \eqref{eq103}  holds for 
all basic tensors $x^* \otimes y^*$. 
Indeed, if $x^*$ or $y^*$ is a biorthogonal functional, then 
obviously 
$\langle B_n, x^* \otimes y^*\rangle  = \sum_{j=1}^n x^* (e_{i_j}^{m_n}) 
y^* (e_{i_j}^{m_n}) =0$ 
for $n$ sufficiently large. 
But also if $\gamma,\beta$ are difference branches of $\B$, then for $n$ 
sufficiently large, the branches will have split at level $m_n$, which implies 
that $\gamma^*(e_{i_j}^{m_n})\beta^* (e_{i_j}^{m_n}) =0$ for all $1\le j\le n$.
Finally, suppose  $\gamma$ is a given branch. 
Then for all $F\in \F$ with $\gamma\in \F$, $F= \{\gamma_1,\ldots,\gamma_n\}$ 
with say $\gamma = \gamma_j$, 
\begin{equation}\label{eq111}
\langle B_F,\gamma^* \otimes \gamma^*\rangle = \gamma_j^*(e_{i_j}^{m_n})
\gamma_j^* (e_{i_j}^{m_n}) =1\ .
\end{equation}
It follows that our cluster point $G_n$ must satisfy
\begin{equation}\label{eq112}
G_2 (\gamma^* \otimes \gamma^*) =1\ .
\end{equation}
Now finally, let $G= G_1 - G_2$. 
Then we have by \eqref{eq101}, \eqref{eq109} and the fact that \eqref{eq103}
holds for $H = G_i$, $i=1,2$, that $G$ satisfies \eqref{eq96} and \eqref{eq97},
and so $G\perp$ $JT_p\Atimes JT_q$, $G(T) =1$, whence $T$ is not nuclear.

Now suppose that $q<p^*$. 
To show that 1 holds, we only need to show that $JT_p$ and $JT_q$ 
satisfy the assumptions of Corollary~7;		
but of course these spaces 
satisfy all the assumptions preceding the final one, so we only need to prove
that 
\begin{equation}\label{eq113} 
\text{Every operator from $(JT_p)^*$ to $(JT_q)^{**}$ has separable range.}
\end{equation}
Let $P:(JT_q)^{**}$ be the projection from $(JT_q)^{**}$ onto $JT_q$, with 
kernel$((JT_q)_*)^\perp$, which we know by Theorem~10	
is isometric to $\ell_\G^q$; for the sake of notational simplicity, let us just
set $((JT_q)_*)^\perp =\ell_\G^q$. 
Now let $T:(JT_p)^* \to (JT_q)^{**}$ be a given operator, and set 
\begin{equation}\label{eq114} 
A = PT \ \text{ and }\ B= (I-P) T\ .
\end{equation}
Obviously $A$ has separable range, and $B$ is an operator from $(JT_p)^*$ 
to $\ell_\G^q$. 
We claim that
\begin{equation}\label{eq115}
B\ \text{ is compact.}
\end{equation}

It suffices to prove that $B^* : \ell_\G^{q^*} \to (JT_p)^{**}$ is compact.
To do that, we shall show that 
\begin{eqnarray}
&&\text{Every normalized weakly null sequence in $(JT_p)^{**}$}\label{eq116}\\
&&\text{has a subsequence equivalent to the $\ell^p$ basis.}\notag
\end{eqnarray} 
This will complete the proof. 
Indeed, to show that $B^*$ is compact, it suffices 
(by reflexivity of $\ell_\G^{q^*}$) 
to show that
\begin{equation}\label{eq117}
\text{given $(x_n)$ a normalized weakly null sequence in $\ell_\G^{q^*}$, 
then $\|B^* (x_n)\| \to0$.}
\end{equation}
Suppose this is not the case.
But if $(x_n)$  fails to satisfy \eqref{eq117}, there exists a subsequence 
$(x'_n)$ such that $(x'_n)$ is equivalent to the $\ell^{q^*}$ basis 
and by \eqref{eq116}, $B^* (x'_n)$ is equivalent to the $\ell^p$ basis.
But our assumption $q<p^*$ is the same as $p<q^*$, and we have thus deduced
that the $\ell^{q^*}$ basis dominates the $\ell^p$ basis, a contradiction.
Of course once we know that $B^*$ and hence $B$, is compact, then $B$ has 
separable range, and thus $T$ has separable range, finishing the proof.

Since $(JT_p)^{**} = JT_p \oplus \ell_\G^p$ and so both factors have the 
desired property, by Theorem~\ref{thm8} part~5, \eqref{eq116} is easily seen,
but for completeness, here's the argument. 
Let $(x_n)$ be a normalized weakly null sequence in $(JT_p)^{**}$. 

We may assume (by passing to a subsequence) that $(x_n)$ is a basic sequence.
Now if there is a subsequence $(x'_n)$ of $(x_n)$ such that $\|Px'_n\|\to0$, 
we easily obtain a further subsequence $(x''_n)$ such that $((I-P)x''_n)$
is equivalent to $(x''_n)$, and hence by passing  to a further 
subsequence $(x''_n)$ if necessary, that $(x'''_n)$ is equivalent to the 
$\ell^p$ basis. 
Thus we may now assume that 
\begin{equation}\label{eq118}  
\text{there is a $\delta>0$ so that $\|Px_n\|  \ge \delta$ for all $n$.}
\end{equation}
But then by Theorem~10		
part 5, we may choose a subsequence 
$(x'_n)$ of $(x_n)$ so that 
\begin{equation}\label{eq119} 
(Px'_n)\ \text{ is equivalent to the $\ell^p$-basis.}
\end{equation}
Now if some subsequence $(x''_n)$ of $(x'_n)$ satisfies that 
$\|(I-P) x''_n\|\to0$, then just as before, we get a subsequence of $(x'_n)$ 
equivalent to the $\ell^p$ basis, so suppose this also is not the case.
But then (after dropping a few terms if necessary), 
there is a $\delta'>0$ such that 
\begin{equation}\label{eq120} 
\|(I-P)x'_n\| \ge \delta'\ \text{ for all }\ n\ ,
\end{equation}
and now by the standard property of the $\ell^p$ basis, there is a 
subsequence $(x''_n)$ such that also 
\begin{equation}\label{eq121}
((I-P) x''_n) \ \text{ is equivalent to the $\ell^1$-basis.}
\end{equation}
But then also by \eqref{eq119}, 
\begin{equation}\label{eq122}
(Px''_n) \ \text{ is equivalent to the $\ell^p$-basis.}
\end{equation}
The continuity of $P$ and $I-P$ and (126),(127) 
now easily yield that $(x''_n)$ 
is equivalent to the $\ell^p$ basis.\qed
\medskip

\begin{remarks} 
1.  Going through the proofs of Theorems~6 and 11  
is a difficult way to see that condition~4 in the Remark (following the 
statement of Theorem~\ref{thm11})  holds if 2 holds, 
for this may more easily be seen directly as follows.
Since $(JT_p)^*$ has the metric approximation property,
\begin{equation}\label{eq123}
JT_p^* \Atimes JT_q^*\ \text{ is isometric to a closed linear subspace of }
I(JT_p,JT_q^*) = (JT_p \vtimes JT_q)^*\ ,
\end{equation}
and thus 
\begin{equation}\label{eq124}
\L_{p,q}\ \text{ is isometric to a quotient space of } (JT_p \vtimes JT_q)^{**},
\end{equation}
(where $\L_{p,q}$ is as in the above remark). 
Using Theorem~10     
part 4 and (77),		
let $Q:(JT_p)^* \to \ell_\G^{p^*}$
be a quotient map and $j:\ell_\G^q\to (JT_q)^{**}$ be an isometric injection. 
Since $p^* \le q$, the $\ell^{**}$-basis dominates the $\ell^q$-basis.
Then given $\eta \in\{ 0,1\}^\G$, there exists a unique linear contraction
$T_\eta :\ell_\G^{p^*} \to \ell_\G^?$ such that 
\begin{equation}\label{eq125}
\text{for all } \alpha <\G,\ (T_\eta) (\alpha) = 1 \text{ if }  
\eta (\alpha) =1,\ (T_\eta) (\alpha)=0\text{ if }\eta (\alpha)=0.
\end{equation}
Now defining $\tilde T_\eta\in \L_{p,q}$ by 
\begin{equation}\label{eq126}
\tilde T_\eta = jT_\eta Q\ ,
\end{equation}
it follows that $\eta \ne\eta \Rightarrow \tilde T_\eta \ne \tilde T_{\eta'}$,
proving that $\card \L_{p,q} = 2^\G$ (and also showing directly that 
$(JT_p \vtimes JT_q)^{**}$ has cardinality $2^\G$ by \eqref{eq124}, thus 
giving that $\ell^1$ embeds in $JT_p\vtimes JT_q$ by \cite{OR}.

2. An inspection of the operator $V$ constructed in the above proof shows 
that its range is actually contained in $(JT_q)_*$, and thus condition
3. of Theorem~\ref{thm10} may be strengthened to 

$3'$. 
{\em There exists an integral operator from $JT_p$ to $(JT_q)$, which is not 
nuclear}\newline

4. It's conceivable that $(JT)^*\Atimes (JT)^*$ (or more generally, 
$(JT_p)^* \otimes (JT_p^*)^*$ for some $1<p,q <\infty$, ($p^*\le q$), 
is a counterexample to Problem~4.
I reduced this to a separable issue, which, however, I cannot decide.
That is, I proved that if $p,q$ are as above and {\em $L^1$ embeds in 
$(JT_p)^* \Atimes (J_q)^*$, then $L^1$ embeds in $(JT_p)_* \Atimes (JT_q)_*$.}
An alternate approach to Problem~3, using embeddings of unconditional families 
in place of embeddings of $L^1$, does not work. 
In fact, it follows that $Y = ((JT)^* \Atimes JT^*)/(JT)_* \Atimes (JT)_*)$ 
is isometric to $\ell_\G^2 \Atimes \ell_\G^2$; the space of ``diagonal'' 
operators in this space is isometric to $\ell_\G^1$. 
Thus, $Y$ contains a subspace isometric to $\ell_\G^1$, whence by the 
lifting property of this space, $(JT)^* \Atimes (JT)^*$ contains a 
subspace isometric to $\ell_\G^1$. 
However $(JT)^*$ has no uncountable unconditional family by 
Theorem~\ref{thm:main}.
\end{remarks}

We conclude  this section with an application of Theorem~\ref{thm8} to the 
relationship of Problem~1 with the RNP. 
We first need 

\begin{lem}\label{lem12} 
$(JT_q)_*$ fails the RNP for all $q$, $1<q<\infty$. 
\end{lem}

\begin{proof} 
We shall show that $(JT_q)_*$ has a ``$\delta$-tree.'' 
We first motivate the construction.
Let $\mu$ be the Cantor measure on the Borel subsets of $\B$ and 
$V:L^1 (\mu) \to (JT_q)^*$ be the operator constructed at the beginning 
of the proof of Theorem~\ref{thm11}. 
Then in fact $V$ is valued in $(JT_q)_*$, and is not representable by a 
Bochner integrable function.
Rather than proving this, however, let us just examine:
$V(\bone)$. Since $(e_j^*)$ is a $w^*$-basis for $(JT_q)^*$, $V(\bone)$
must have an expansion; $V(\bone) = \sum_{j=1}^\infty c_j e_j^*$, the series
converging in the $w^*$-topology.
The coefficients $(c_j)$ are determined by: $c_j = \langle V(\bone),e_j\rangle$
for all $j$. 
Let us instead, label the coefficient corresponding to $e_j$ as $c_\alpha$, 
where $j=\tau(\alpha)$. 
Then for any $\alpha \in \D$, 
\begin{equation}\label{eq127}
\langle V\bone,e_{\tau(\alpha)} \rangle 
= \int_\B w^* (e_{\tau(\alpha)}) \,d\mu(w) = \frac1{2^{|\alpha|}}\ ,
\end{equation}
because $w^* (e_{\tau(\alpha)}) =1$ if $w\in \B_\alpha$ and 0 otherwise.

But we have that 
\begin{equation}\label{eq128} 
\sum_{\alpha\in\D} \frac1{2^{|\alpha|}} e_{\tau (\alpha)}^*\ \text{ converges 
in norm.}
\end{equation}
Indeed, let $y_n$ be defined by 
\begin{equation}\label{eq129}
y_n = \sum_{j=1}^{2^n} \frac1{2^n} (e_j^n)^* = \frac1{2^n} 
\sum_{|\alpha|=n} e_{\tau(\alpha)}^*\ .
\end{equation}
Then 
\begin{equation}\label{eq130}
\|y_n\| = \frac1{2^{n/q^*}} \ 
\qquad \text{ by \eqref{eq84}.}
\end{equation}
Thus
\begin{equation}\label{eq131}
\sum_{n=0}^\infty \|y_n\| = \frac1{1-\frac1{2^{n/q^*}}} 
= \frac{2^{n/q^*}}{2^{n/q^*}-1} \defeq c_q <\infty\ ,
\end{equation}
and of course $\sum_{n=0}^\infty y_n = \sum_{\alpha\in\infty} 
\frac1{2^{|\alpha|}} e_{q(\alpha)}^*$.

It is then easily verified that for any $\alpha$,
\begin{equation}\label{eq132}
T(\chi_{\B_\alpha}) = \sum_{\gamma \in \B_\alpha} \frac1{|\alpha|} 
e_{\tau(\alpha)}^*
\end{equation}
and this series converges in norm, by exactly the same reason we gave 
for \eqref{eq128}.

Now define $(t_\alpha)_{\alpha\in\D}$ by 
\begin{equation}\label{eq133}
t_\alpha = 2^{|\alpha|} \sum_{\gamma\in \B_\alpha} \frac1{2^{|\alpha|}}
e_{\tau(\gamma)}^* = T(2^{|\alpha|} \chi_{\B_\alpha} \ .
\end{equation}
Since $[e_{\tau(\gamma)}^* :\gamma \in \B_\alpha]$ is canonically isometric to 
$(JT_q)_*$, we thus have that 
\begin{equation}\label{eq134}  
\|t_\alpha\| \le c_q\ \text{ for all }\ \alpha\ .
\end{equation}
It follows by \eqref{eq133} and the linearity of $T$ (or by direct 
computation) that for all $\alpha \in \D$, 
\begin{equation}\label{eq135} 
t_\alpha = \frac12 (t_{\alpha_0} + t_{\alpha_1})
\end{equation}
and 
\begin{equation}\label{eq136}
\|t_{\alpha_0} - t_{\alpha_1} \| \ge \|e_{\alpha_0}^* - e_{\alpha_1}^*\| 
= 2^{1/q^*}\ .
\end{equation}
Thus $(t_\alpha)_{\alpha\in\D} $ is a bounded ``$\delta$-tree'' with 
$\delta = 2^{1/q^*}$, and so $(JT_q)_*$ fails the RNP
(cf. \cite{DU}).
\end{proof}

\begin{cor}\label{cor13} 
Let $1<p<\infty$ and let $Y_p = (JT_p)_* \oplus JT_p$.
\begin{itemize}
\item[1.] $Y_p$ and $Y_p^*$ fail the RNP.
\item[2.] $\ell^1$ does not embed in $K(Y_p)$ if $1<p<2$.
\item[3.] $\ell^1$ embeds in $K(Y_p)$ and moreover there exists an integral
non-nuclear operator on $Y_p^*$ if $2\le p<\infty$.
\end{itemize}
\end{cor}

\begin{proof}
1. $Y_p$ fails the RNP by the preceding result and 
\begin{equation}\label{eq137}
Y_p^* = JT_p \oplus (JT_p)^*
\end{equation}
fails the RNP because $(JT_p)^*$ does.

Next, we have for any $p$ that 
\begin{equation}\label{eq138}
K(Y_p) = 
\left[ (JT_p)_* \vtimes JT_p\right] \oplus 
\left[ (JT_p)^* \otimes JT_p\right] \oplus 
\left[ (JT_p)_* \vtimes (JT_p)^*\right] \oplus 
\left[ JT_p \vtimes JT_p\right]\ .
\end{equation} 
$\ell^1$ does not embed in the first three summands in \eqref{eq138} 
by Corollary~\ref{cor:Lewis},  
because in each case, the dual  of one of the factors in the injective 
tensor product has the RNP (and of course the factors all have the map,
and their duals do not contain $\ell^1$ isomorphically).

If $1<p<2$, then $p^* >p$, and so $\ell^1$ does not embed in the 
fourth summand by Theorem~\ref{thm11}, so $\ell^1$ does not embed in $K(Y_p)$.

If $2\le p$, then $p^* \le p$, and thus $\ell^1$ embeds in the fourth 
summand by Theorem~\ref{thm11}, part~3 of which also shows part~3 of the 
Corollary, using \eqref{eq138}.
\end{proof}

\begin{Rem}
Notice by Lemma~12      
that we also obtain an integral operator
from $JT_p$ into $(JT_p)_*$ which is not nuclear 
if $2\le p<\infty$, while if $1<p<2$, every integral operator 
from $JT_p$ into $(JT_p)^*$ is nuclear by Theorem~\ref{thm11}. 
\end{Rem}

\appendix{}
\section*{Appendix}

We give a somewhat new proof of \eqref{eq8} via 

\begin{prop}\label{prop:A1}
Suppose that $Z$ is a $\L_\infty$  space which is isomorphic to a quotient
of a subspace of a Banach space $X$. 
Then $Z^*$ is isomorphic to a subspace of $X^*$.
\end{prop}

Thus to obtain \eqref{eq8}, if $\ell^1$ is isomorphic to a subspace of $X$, 
$C([0,1])$ is isomorphic to a quotient of that subspace, and hence the 
Proposition applies.
(For properties of $\L_\infty$ spaces, see \cite{LP} and \cite{LR}; the 
definition will appear in our proof.) 

\begin{proof}[Proof of Proposition~\ref{prop:A1}]
Let $\tilde X$ be a (closed linear) subspace of $X$ and $T :\tilde X\to Z$ 
a surjective bounded linear map, with $\|T\|=1$. 
We may choose a ``$\P_1$'' space $W$ with $Z\subset W$ (e.g., $W=\ell^\infty 
(\Ba X^*))$, and so we can then choose $\tilde T:X\to W$ with 
$\|\tilde T\|=1$ such that $\tilde T$ extends $T$.
Now choose $\lambda >1$ such that for all finite dimensional subspaces 
$E$ of $Z$, there exists a subspace $F_E\supset E$ with 
\begin{equation}\label{eq139} 
d(F_E,\ell_n^\infty) \le \lambda\ ,\ \text{ where }\ \dim F_E<\infty 
\end{equation}
(and the first term in \eqref{eq139} denotes the Banach-Mazur distance from 
$F_E$ to $\ell_n^\infty$). 
Therefore we may choose $P_E :W\to Z$ with 
\begin{equation}\label{eq140} 
\|P_E\|\le \lambda\ ,\quad P_E^2 = P_E\ ,\ \text{ and }\ P|_{F_E} = I|_{F_E}\ .
\end{equation}

Now let $\D$ be the family of all finite dimensional subspaces of $Z$, 
directed by reverse inclusion.
Then a compactness argument, using the fact that $\lambda \Ba Z^*$ is 
$w^*$-compact, shows there is a linear operator $S:Z^* \to W^*$ so that 
\begin{equation}\label{eq141} 
\|S\| \le \lambda\ \text{ and for all }\  z\in Z\ \text{ and }\ 
z^* \in Z^*\ ,\ \langle Sz^*,z\rangle = \langle z^*,z\rangle\ .
\end{equation}

Thus if $R:W^*\to Z^*$ denotes the canonical restriction map, we have by 
\eqref{eq141} that 
\begin{equation}\label{eq142} 
RS = I_{Z^*}\ .
\end{equation}
Then setting $Y = S(Z^*)$, $Y$ is $\lambda$-isomorphic to $Z^*$, 
since \eqref{eq142} gives that $S$ is an isomorphism with the inverse 
of $S:Z^*\to Y$ given by ${R|Y }$. 
Moreover, 
\begin{equation}\label{eq143}
\begin{array}{l}
\text{Given $\ep>0$ and $y\in Y$, there exists a $z\in Z$ with $\|z\|=1$}\\
\text{such that $|y(z)| > \frac1{\lambda} \|y\| -\ep$.}
\end{array}
\end{equation}
(It also follows from \eqref{eq142} that $Y$ is complemented in $Z^*$, by 
$Z^\perp$, but we don't use this.)

Now \eqref{eq143}  implies that 
\begin{equation}\label{eq144} 
\tilde T^* \mid Y \text{ is an isomorphism, with $\tilde T(Y)$ being 
$\lambda\| (T^*)^{-1}\|$-isomorphic to $Z^*$.}
\end{equation}
Indeed, given $y\in Y$ and $\ep>0$, choose $z\in Z$ satisfying \eqref{eq143},
then choose $x\in\tilde X$ with $\|x\| \le \|(T^*)^{-1}\| +\ep$ so that 
$Tx = z$. 
Thus letting $\tau = \frac1{\|(T^*)^{-1}\|} +\ep$, 
\begin{equation}\label{eq145} 
\begin{split}
\|(\tilde T^*) (y) \| 
& = \tau |\langle \tilde T^* y,x\rangle |\\
& = \tau |\langle y, \tilde Tx\rangle | = \tau |\langle y,Tx\rangle|\\
& = \tau |\langle y,z\rangle|\\
&\ge \frac{\tau}{\lambda} \|y\| - \ep\ .
\end{split}
\end{equation}
Since $\ep>0$ was arbitrary, \eqref{eq144} is proved.\qed

\begin{Rem}
We have kept track of the constants in the above proof because it then 
yields the following result:
{\em Suppose $X$ and $Z$ are Banach spaces so that $Z$ is an ``$L^1$-predual,''
meaning that $Z^*$ is isometric to $L^1(\mu)$ for some (not necessarily 
$\sigma$-finite) measure $\mu$, and suppose for all $\ep>0$, there exists 
a subspace $\tilde X$ of $X$ so that $Z$ is $1+\ep$-isomorphic to a 
quotient space of $\tilde X$. 
Then $Z^*$ is $1+\ep$-isomorphic  to a subspace of $X^*$ for all $\ep>0$.}
\end{Rem}

Indeed, we need only apply our proof, using the standard fact (as follows 
from local reflexivity) that $Z$ is an $\L_{\infty,\lambda}$ space for all 
$\lambda>1$. 
Now assuming $\ell^1$ embeds in $X$, then by a result of James \cite{J1}, given 
$\ep>0$, there is a subspace $\tilde X$ of $X$ which $(1+\ep)$-isomorphic 
to $\ell^1$, and consequently $C([0,1])$ is $(1+\ep)$-isomorphic to a quotient
space of $\tilde X$. 
Hence we deduce that {\em for all $\ep>0$, $(C(\Lambda))^*$ is 
$(1+\ep)$-isomorphic to a subspace of $X^*$.}

As we show in Proposition~17,     
James' result that $\ell^1$ is not 
distortable, holds for the spaces $\ell^1_\kappa$ as well, $\kappa$ an 
infinite cardinal.
We thus obtain the following generalization of the above quantitative 
version of \eqref{eq8}. 

\begin{thm}\label{thm15}
Let $\kappa$ be an infinite cardinal number, and suppose $X$ contains a 
subspace isomorphic  to $\ell^1_\kappa$. 
Then given $Z$ an $L^1$ predual of density character at most $\kappa$, 
$X^*$ contains a subspace $(1+\ep)$-isomorphic  to $Z^*$ for all $\ep>0$.
\end{thm}

The following result, which reduces to Pe{\l}czy\'nski's theorem for 
$\kappa = \aleph_0$, is an immediate consequence.

\begin{cor}\label{cor16}
Let $\kappa$ be an infinite cardinal number and suppose $X$ contains a 
subspace isomorphic to $\ell^1_\kappa$. 
Then for all $\ep>0$, $X^*$ contains a subspace $(1+\ep)$-isomorphic to 
$[C(\{0,1\}^\kappa )]^*$.
\end{cor}

\end{proof}

Finally, we give the extension of James' theorem to the spaces $\ell^1_\kappa$
(which is apparently a new result).

\begin{prop}\label{prop:James-extension}
Let $\kappa$ be an infinite cardinal number, and assume $X$ contains a 
subspace isomorphic to $\ell^1_\kappa$. 
Then $X$ contains a subspace $(1+\ep)$-isomorphic to $\ell^1_\kappa$ for 
all $\ep>0$.
\end{prop}

\begin{proof}
We identify cardinals with initial ordinals. 
Let then $[e_\alpha)_{\alpha <\kappa}$ be a normalized basis of cardinality 
$\kappa$ in $X$, equivalent to the $\ell^1_\kappa$ basis.
For each cardinal $\alpha <\kappa$, define $\delta_\alpha$ by 
\begin{equation}\label{eq151}
\begin{cases}
\ds\delta_\alpha =\sup\{\delta >0:\|\sum_{\gamma\ge \alpha} c_\gamma e_\gamma\|
\ge \delta\ \text{ for}\\
\ds \text{all families of scalars }\ (c_\gamma)_{\gamma\ge \alpha}\
\text{ such that }\  \sum_{\gamma\ge\alpha} |c_\gamma| = 1\}
\end{cases}
\end{equation}
It is obvious that $\alpha\to \delta_\alpha$ is an increasing function.
Hence 
\begin{equation}\label{eq152}
\lim_{\alpha\to\kappa} \delta_\alpha  \defeq \delta\ \text{ exists.}
\end{equation}
(Of course if $\kappa$ is of uncountable cofinality, then 
$\delta_\alpha=\delta$ for all $\alpha$ sufficiently large; however  this fact
is irrelevant for the proof.) 

Now let $0<\eta<\delta$ be given. 
It follows easily by induction and the fact that the family of all finite 
subsets of $\kappa$ also has cardinality $\kappa$, and also because 
$\card \gamma < \kappa$ for $\gamma <\kappa$ and hence $\card \{\alpha :\gamma
\le \alpha < \kappa\}=\kappa$ for all $\gamma <\kappa\}$, that we 
may choose a family $(f_\alpha)_{\alpha <\kappa}$ of finite linear 
combinations of the $e_\alpha$'s such that for each $\alpha$,
there exist ordinals $\alpha \le a_\alpha \le b_\alpha$ and $c_\gamma$'s, 
only finitely many non-zero, with 
\begin{gather}
f_\alpha = \sum_{a_\alpha \le\gamma\le b_\alpha} c_\gamma e_\gamma
\tag{153i}\label{tag153i}\\
\noalign{\vskip6pt}
\|f_\alpha\| =1\quad\text{and}\quad \sum_{a_\alpha \le\gamma\le b_\alpha} 
|c_\gamma| > \frac1{\delta_0+\eta}\tag{153ii}\label{tag153ii}\\
\noalign{\vskip6pt}
\delta_{\alpha_0} > \delta -\eta\tag{153iii}\label{tag153iii}\\
\noalign{\vskip6pt}
\text{for all }\ \alpha <\alpha' <\kappa\ ,\quad b_\alpha < a'_\alpha\ .
\tag{153iv}\label{tag153iv}
\end{gather}
Now let $(x_\alpha)_{\alpha<\kappa}$ be a family of scalars, only finitely 
many non-zero, such that $\sum_{\alpha <\kappa} |x_\alpha| =1$. 
Thus we have 
\addtocounter{equation}{1}
\begin{equation}
\begin{split}
\|\sum_\alpha x_\alpha f_\alpha\| 
& = \| \sum_\alpha x_\alpha  \sum_{a_\alpha \le\gamma\le b_\alpha} 
c_\gamma e_\gamma\|\\
\noalign{\vskip6pt}
& \ge \delta_{\alpha_0}  \sum_\alpha |x_\alpha| 
\sum_{a_\alpha\le\gamma\le b_\alpha}|c_\gamma| \ \text{ by definition of }\ 
\delta_{\alpha_0}\\
\noalign{\vskip6pt}
&> \frac{\delta_{\alpha_0}}{\delta_0 +\eta} \sum |x_\alpha|
\ \text{ by \eqref{tag153ii}} \\
\noalign{\vskip6pt}
& > \frac{\delta-\eta}{\delta +\eta} \ \text{ by \eqref{eq152}.}
\end{split}
\end{equation}
Thus given $\ep >0$, we may choose $\eta$ so small that 
$\frac{\delta-\eta}{\delta+\eta} > \frac1{1+\ep}$, proving 
Proposition~\ref{prop:James-extension}.
\end{proof}

\begin{Rem}
Similar reasoning shows that the ``predual'' formulation of 
Proposition~17 also holds.
That is, 
\end{Rem}

\begin{prop}\label{prop:predual}
Let $\kappa$ be a infinite cardinal number, and assume $X$ contains 
a subspace isomorphic to $c_0(\kappa)$.
Then $X$ contains a subspace $(1+\ep)$-isomorphic to 
$c_0(\kappa)$ for all $\ep>0$.
\end{prop}


\end{document}